\documentclass[11pt]{article}
\usepackage{srcltx}
\usepackage{epsfig}
\usepackage{url}
\def\Fro{{\hbox{\scriptsize\rm Fro}}}
\def\cL{{\cal L}}

\def\nuc{{\hbox{\scriptsize\rm nuc}}}

\def\cC{{\cal C}}

\def\VI{{\hbox{\rm VI}}}
\def\epsilonvi{\epsilon_{\hbox{\scriptsize\rm vi}}}
\def\epsilonsad{\epsilon_{\hbox{\scriptsize\rm sad}}}
\def\Res{{\hbox{\rm Res}}}
\def\Gap{{\hbox{\rm Gap}}}

\usepackage{psfrag}
\usepackage{graphicx}
\usepackage[latin1]{inputenc}
\usepackage{subfigure}
\usepackage{hyperref}
\usepackage{mathtools}
\usepackage{amsfonts,amsmath,amssymb,amsthm}
\usepackage{color}
\usepackage{booktabs}

\usepackage{float,enumitem}
\setlist{itemsep=3pt,parsep=0pt,topsep=3pt}

\floatstyle{ruled}
\newfloat{algorithm}{tbp}{loa}
\floatname{algorithm}{Algorithm}

\DeclareMathOperator{\rank}{rank}

\DeclareMathOperator*{\argmax}{Arg\,max}
\DeclareMathOperator*{\argmin}{Arg\,min}

\newtheorem{thm}{Theorem}

\newtheorem{lemma}[thm]{Lemma}

\DeclareMathOperator{\Prox}{Prox}

\def\qed{\ \hfill$\square$\par\smallskip}

\def\Tr{{\mathop{\hbox{\rm Tr}}}}
\def\Opt{{\mathop{\hbox{Opt}}}}

\def\rank{{\mathop{\hbox{\rm Rank}}}}

\def\Diag{{\hbox{\rm Diag}}}

\newcommand{\cN}{I\!\! N}

\newcommand{\half}{ \mbox{\small$\frac{1}{2}$}}

\newcommand{\be}{\begin{eqnarray}}
\newcommand{\ee}[1]{\label{#1}\end{eqnarray}}

\newcommand{\ese}{\end{eqnarray*}}
\newcommand{\bse}{\begin{eqnarray*}}
\newcommand{\rf}[1]{~(\ref{#1})}

\newtheorem{proposition}{Proposition}
\newtheorem{corollary}{Corollary}
\newtheorem{theorem}{Theorem}
\newtheorem{remark}{Remark}
\def\argmin{\mathop{\hbox{\rm argmin$\,$}}}
\def\Argmin{\mathop{\hbox{\rm Argmin$\,$}}}
\def\argmax{\mathop{\hbox{\rm argmax$\,$}}}

\def\cA{{\cal A}}

\def\cN{{\cal N}}

\def\bR{{\mathbf{R}}}

\definecolor{darkmagenta}{RGB}{125,38,205}

\begin{document}
\title{Solving Variational Inequalities with Monotone Operators on Domains Given by Linear Minimization Oracles}
\author{Anatoli Juditsky\thanks{LJK, Universit\'e Grenoble Alpes, B.P. 53, 38041 Grenoble Cedex 9, France, {\tt Anatoli.Juditsky@imag.fr}}
\and Arkadi Nemirovski\thanks{Georgia Institute
 of Technology, Atlanta, Georgia
30332, USA, {\tt nemirovs@isye.gatech.edu}\newline
Research of the first author was supported by the CNRS-Mastodons project GARGANTUA,
and the LabEx PERSYVAL-Lab (ANR-11-LABX-0025). Research of
the second author was supported by the NSF grants CMMI 1232623 and CCF 1415498.}
}
\maketitle
\begin{abstract}
 The standard algorithms for solving large-scale convex-concave saddle point problems, or, more generally, variational inequalities with monotone operators, are {\sl proximal type} algorithms which at every iteration need to compute a {\sl prox-mapping}, that is, to minimize over problem's domain $X$ the sum of a linear form and the specific convex {\sl distance-generating function} underlying the algorithms in question. (Relative) computational simplicity of prox-mappings, which is the standard requirement when implementing proximal algorithms, clearly implies the possibility to equip $X$ with a relatively computationally cheap {\sl Linear Minimization Oracle} (LMO) able to minimize over $X$ linear forms. There are, however, important situations where a cheap LMO indeed is available, but where no proximal setup with easy-to-compute prox-mappings is known.  This fact motivates our goal in this paper, which is to develop techniques for solving variational inequalities with monotone operators on domains given by Linear Minimization Oracles. The techniques we are developing can be viewed as a substantial extension of the proposed in \cite{DualSubgr} method of nonsmooth convex minimization over an LMO-represented domain.\end{abstract}
 \section{Introduction}\label{sec:intro}
 The majority of First Order methods (FOM's) for large-scale convex minimization (and {\sl all} known to us FOM's for large-scale convex-concave saddle point problems and variational inequalities with monotone operators) are of {\sl proximal} type: at a step of the algorithm, one needs to compute {\sl prox-mapping} -- to minimize over problem's domain the sum of a linear function and  a specific for the algorithm strongly convex {\sl distance generating} function (d.-g.f.), in the simplest case, just squared Euclidean norm. As a result,  the  practical scope of proximal algorithms is restricted to {\sl proximal-friendly} domains -- those allowing for d.-g.f.'s with not too expensive computationally  prox-mappings.
 What follows is motivated by the desire to develop FOM's for solving convex-concave saddle point problems on bounded domains with ``difficult geometry'' -- those for which no d.-g.f.'s resulting in nonexpensive prox-mappings (and thus no ``implementable'' proximal methods) are known.  In what follows, we relax the assumption on problem's domain to be proximal-friendly to the weaker assumption to admit computationally nonexpensive {\sl Linear Minimization Oracle} (LMO) -- a routine capable to minimize a linear function over the domain. This indeed is a relaxation: to minimize within a desired, whatever high, accuracy a linear form over a bounded proximal-friendly domain is the same as to minimize over the domain the sum of large multiple of the form and the d.-g.f. Thus, proximal friendliness implies existence of a nonexpensive LMO, but not vice versa. For example, when the domain is the ball  $B_n$ of nuclear norm in $\bR^{n\times n}$, computing prox-mapping, for all known proximal setups, requires full singular value decomposition of an $n\times n$ matrix, which can be prohibitively time consuming when $n$ is large. In contrast to this, minimizing a linear form over $B_n$ only requires finding the leading singular vectors of an $n\times n$ matrix, which is much easier than full-fledged singular value decomposition.
 \par
Recently, there was significant interest in solving convex minimization problems on  domains given by LMO's. The emphasis in this line of research
is on smooth/smooth norm-regularized convex minimization \cite{DHM.7,Fre,HDPDM.11,Jag1.16,Jag2.17,SS.32}, where the main ``working horse'' is the classical Conditional Gradient (a.k.a. Frank-Wolfe) algorithm originating from \cite{FW} and intensively studied in 1970's (see  \cite{DemRub,Dunn,PsheD} and references therein). Essentially, Conditional Gradient is the only traditional convex optimization technique capable to handle convex minimization problems on LMO-represented domains. In its standard form,  Conditional Gradient algorithm, to the best of our knowledge, is not applicable beyond the smooth minimization setting; we are not aware of any attempt to apply this algorithm even to the simplest -- bilinear -- saddle point problems. The approach proposed in this paper is different and is inspired by our recent paper \cite{DualSubgr}, where a method for {\sl nonsmooth} convex minimization over an LMO-represented convex domain was developed. The latter method utilizes  {\sl Fenchel-type} representations of the objective in order to pass from the problem of interest to its special dual. In many important cases the domain of the dual problem is proximal-friendly, so that the dual problem can be solved by proximal FOM's. We then use the machinery of {\sl accuracy certificates} originating from \cite{NOR} allowing to recover a good solution to the problem of interest from the information accumulated when solving the dual problem. In this paper we follow the same strategy in the context of variational inequalities (v.i.'s) with monotone operators (this covers, in particular, convex-concave saddle point problems). Specifically,   we introduce the notion of a {\sl Fenchel-type  representation} of a monotone operator, allowing to associate with the v.i. of interest its dual, which is again a v.i. with monotone operator with the values readily given by the representation {\sl and the LMO representing the domain of the original v.i.} Then we solve the dual v.i. (e.g., by a proximal-type algorithm) and use the machinery of accuracy certificates to recover a good solution to the v.i. of interest from the information gathered when solving the dual v.i.
\par
The main body of the paper is organized as follows. Section \ref{secprel} outlines the background of convex-concave saddle point problems, variational inequalities with monotone operators and accuracy certificates. In section \ref{secrepr}, we introduce the notion of a Fenchel-type representation of a monotone operator and the induced by this notion concept of v.i. dual to a given v.i. This section also contains a simple fully algorithmic ``calculus'' of Fenchel-type representations of monotone operators: it turns out that basic monotonicity-preserving
operations with these operators (summation, affine substitution of argument, etc.) as applied to operands given by Fenchel-type representations yield similar representation for the result of the operation. As a consequence, our abilities to operate numerically with Fenchel-type representations of monotone operators are comparable with our abilities to evaluate the operators themselves. Section \ref{secmain} contains our main result -- Theorem \ref{propmain}. It shows how information collected when solving the dual v.i. to some accuracy, can be used to build an approximate solution of the same accuracy to the primal v.i. In Section \ref{secmain} we present a self-contained description of two well known proximal type algorithms for v.i.'s with monotone operators -- Mirror Descent (MD) and Mirror Prox (MP) -- which indeed are capable to collect the required information. Section \ref{secmodif} is devoted to some modifications of our approach as applied to an affine monotone operator. In the  concluding Section \ref{secill}, we illustrate the proposed approach by applying it to the ``matrix completion problem with spectral norm fit'' -- to the problem
\[
\min_{{u\in\bR^{n\times n},\atop\|u\|_\nuc\leq1}}\|\cA u -b\|_{2,2},
\] where $\|x\|_\nuc=\sum_i\sigma_i(x)$ is the nuclear norm, $\sigma(x)$ being the singular spectrum of $x$, $\|x\|_{2,2}=\max_i\sigma_i(x)$ is the spectral norm, and $u\mapsto\cA u$
is a linear mapping from  $\bR^{n\times n}$ to $\bR^{m\times m}$.

\section{Preliminaries}\label{secprel}
\paragraph{Variational inequalities and related accuracy measures.} Let $Y$ be a nonempty closed convex set in Euclidean space $E_y$  and $H(y):Y\to E_y$ be a monotone operator:
$$
\langle H(y)-H(y'),y-y'\rangle \geq 0\,\,\forall y,y'\in Y.
$$
The variational inequality (v.i.) associated with $(H,Y)$ is
$$
\hbox{find $y_*\in Y$}: \langle H(z),z-y_*\rangle \geq 0\,\,\forall z\in Y;\eqno{\VI(H,Y)}
$$
(every) $y_*\in Y$ satisfying the target relation in $\VI(H,Y)$ is called a weak solution to the v.i.; when $Y$ is convex and {\sl compact}, and $H$ is monotone on $Y$, weak solutions always exist. A strong solution  to v.i. is a point $y_*\in Y$ such that   $\langle H(y_*),y-y_*\rangle \geq0$ for all $y\in Y$; from the monotonicity of $H$ is follows that a strong solution is a weak one as well.  Note that when $H$ is monotone {\sl and continuous} on $Y$ (this is the only case we will be interested in), {\sl weak solutions are exactly the  strong solutions}.\par
The accuracy measure naturally quantifying the inaccuracy of a candidate solution $y\in Y$ to $\VI(H,Y)$ is the {\sl dual gap function}
$$
\epsilonvi(y|H,Y)=\sup_{z\in Y}\langle H(z),y-z\rangle;
$$
this (clearly nonnegative for $y\in Y$) quantity is zero if and only if $y$ is a weak solution to the v.i.\par
We will be interested also in the {\sl Special case} where $Y=V\times W$ is the direct product of nonempty convex {\sl compact} subsets $V\subset E_v$ and $W\subset E_w$ of Euclidean spaces $E_v$, $E_w$, and $H$ is associated with Lipschitz continuous function $f(v,w):Y=V\times W\to\bR$ convex in $v\in V$ and concave in $w\in W$:
$$
H(y=[v;w])=[H_v(v,w);H_w(v,w)] \hbox{\ with\ } H_v(v,w)\in\partial_v f(v,w),\,\,H_w(v,w)\in\partial_w[-f(v,w)].
$$
We can associate with the Special case two optimization problems
$$
\begin{array}{rcll}
\Opt(P)&=&\min_{v\in V}\left[\overline{f}(v)=\sup_{w\in W} f(v,w)\right]&(P)\\
\Opt(D)&=&\max_{w\in W}\left[\underline{f}(w)=\inf_{v\in V} f(v,w)\right]&(D)\\
\end{array};
$$
under our assumptions ($V,W$ are convex and compact, $f$ is continuous convex-concave) these problems are solvable with equal optimal values. We associate with a pair $(v,w)\in V\times W$ the {\sl saddle point inaccuracy}
$$
\epsilonsad(v,w|f,V,W)=\overline{f}(v)-\underline{f}(w)=[\overline{f}(v)-\Opt(P)]+[\Opt(D)-\underline{f}(w)].
$$

\paragraph{Accuracy certificates.} Given $Y,\;H$, let us call a collection $\cC^N=\{y_t\in Y,\,\lambda_t\geq0,\,H(y_t)\}_{t=1}^N$ with $\sum_t\lambda_t=1$, an {\sl $N$-step accuracy certificate}. For $Z\subset Y$, we call the quantity
$$
\Res\left(\cC^N|Z\right)=\sup_{y\in Z}\sum_{t=1}^N\lambda_t\langle H(y_t),y_t-y\rangle
$$
the {\sl resolution} of the certificate $\cC^N$ w.r.t. $Z$.\par
Let us make two observations coming back to \cite{NOR}:
\begin{lemma}\label{lem1} Let $Y$ be a closed convex set in Euclidean space $E_y$, $H$ be a monotone operator on $Y$, and $\cC^N=\{y_t\in Y,\lambda_i\geq0,H(y_t)\}_{t=1}^N$ be an accuracy certificate. Setting
$$
\widehat{y}=\sum_{t=1}^N\lambda_ty_t,
$$
we have $\widehat{y}\in Y$, and for every nonempty closed convex subset $Y'$ of $Y$ it holds
\begin{equation}\label{eq1}
\epsilonvi(\widehat{y}|H,Y')\leq \Res\left(\cC^N|Y'\right)
\end{equation}
In the Special case we have also
\begin{equation}\label{eq2}
\epsilonsad(\widehat{y}|f,V,W)\leq \Res\left(\cC^N|Y\right)
\end{equation}
\end{lemma}
{\bf Proof.} For $z\in Y'$ we have
$$
\begin{array}{rl}
\langle H(z),z-\sum_t\lambda_ty_t\rangle =&\sum_t\lambda_t\langle H(z),z-y_t\rangle\\
&\hbox{[since $\sum_t\lambda_t=1$]}\\
\geq&\sum_t\lambda_t\langle H(y_t),z-y_t\rangle \\
&\hbox{[since $H$ is monotone and $\lambda_t\geq0$]}\\
\geq&-\Res\left(\cC^N|Y'\right)\\
&\hbox{[by definition of resolution]}\\
\end{array}
$$
Thus, $\langle H(z),\widehat{y}-z\rangle \leq \Res\left(\cC^N|Y\right)$ for all $z\in Y'$, and (\ref{eq1}) follows. In the Special case, setting $y_t=[v_t;w_t]$, $\widehat{y}=[\widehat{v};\widehat{w}]$, for every $y=[v;w]\in Y=V\times W$ we have
$$
\begin{array}{rcl}
\Res(\cC^N|Y)& \geq& \sum_t\lambda_t\langle H(y_t),y_t-y\rangle\\
&&\hbox{[by definition of resolution]}\\
&=&\sum_t\lambda_t\left[\langle H_v(v_t,w_t),v_t-v\rangle +\langle H_w(v_t,w_t),w_t-w\rangle\right]\\
&\geq&\sum_t\lambda_t\left[[f(v_t,w_t)-f(v,w_t)]+[f(v_t,w)-f(v_t,w_t)]\right]\\
&&\hbox{[by origin of $H$ and since $f(v,w)$ is convex in $v$ and concave in $w$]}\\
&=&\sum_t\lambda_t[f(v_t,w)-f(v,w_t)]\\
&\geq& f(\widehat{v},w)-f(v,\widehat{w})\\
&&\hbox{[since $f(v,w)$ is convex in $v$ and concave in $w$]}
\end{array}
$$
Since the resulting inequality holds true for all $v\in V$, $w\in W$, we get $\overline{f}(\widehat{v})-\underline{f}(\widehat{w})\leq \Res(\cC^N|Y)$, and (\ref{eq2}) follows. \qed
Lemma \ref{lem1} can be partially inverted in the case of {\sl skew-symmetric} operator $H$, that is,
\begin{equation}\label{skew}
H(y)=a+Sy
\end{equation}
with {\sl skew-symmetric} $(S=-S^*)$ \footnote{From now on, for a linear mapping $x\mapsto Bx: E\to F$, where $E,F$ are Euclidean spaces, $B^*$ denotes the conjugate of $B$, that is, a linear mapping $y\mapsto B^*y:F\to E$ uniquely defined by the identity $\langle Bx,y\rangle = \langle x,B^*y\rangle $ for all $x\in E,$ $y\in F$.} linear operator $S$. A skew-symmetric $H$ clearly satisfies the identity
$$
\langle H(y),y-y'\rangle =\langle H(y'),y-y'\rangle, \,\,y,y'\in E_y.
$$
\begin{lemma}\label{lem2} Let $Y$ be a convex compact set in Euclidean space $E_y$, $H(y)=a+Sy$ be skew-symmetric, and let $\cC^N=\{y_t\in Y,\lambda_t\geq0,H(y_t)\}_{t=1}^N$
be an accuracy certificate. Then for $\widehat{y}=\sum_t\lambda_t y_t$ it holds
\begin{equation}\label{eq3}
\epsilonvi(\widehat{y}|H,Y)=\Res\left(\cC^N|Y\right).
\end{equation}
\end{lemma}
{\bf Proof.} We already know that $\epsilonvi(\widehat{y}|H,Y)\leq\Res\left(\cC^N|Y\right)$. To prove the inverse inequality, note that for every $y\in Y$ we have
\bse
&&\begin{array}{ll}
\epsilonvi(\widehat{y}|H,Y)
\geq\langle H(y),\widehat{y}-y\rangle
=&
\langle H(\widehat{y}),\widehat{y}-y\rangle\\
&\mbox{[since $H$ is skew-symmetric]}
\end{array}\\
&&\begin{array}{rll}
=\langle a,\widehat{y}-y\rangle - \langle S\widehat{y},y\rangle +\langle S\widehat{y},\widehat{y}\rangle
=&\langle a,\widehat{y}-y\rangle - \langle S\widehat{y},y\rangle
=&\sum_t\lambda_t[\langle a,y_t-y\rangle -\langle Sy_t,y\rangle]\\
&\mbox{[due to $S^*=-S$]}&\mbox{[due to $\widehat{y}=\sum_t\lambda_ty_t$ and $\sum_t\lambda_t=1$]}
\end{array}
\\
&&\begin{array}{rl}
=&\sum_t\lambda_t[\langle a,y_t-y\rangle +\langle Sy_t,y_t-y\rangle]
=\sum_t\lambda_t\langle H(y_t),y_t-y\rangle.\\
&\mbox{[due to $S^*=-S$]}
\end{array}
\ese
Thus, $\sum_t\lambda_t\langle H(y_t),y_t-y\rangle\leq\epsilonvi(\widehat{y}|H,Y)$ for all $y\in Y$, so that $\Res\left(\cC^N|Y\right)\leq \epsilonvi(\widehat{y}|H,Y)$. \qed
\begin{corollary}\label{newcor}  Assume we are in the Special case, so that $Y=V\times W$ is a direct product of two convex compact sets, and the monotone operator $H$ is associated with a convex-concave function $f(v,w)$. Assume also that $f$ is bilinear: $f(v,w)=\langle a,v\rangle +
 \langle b,w\rangle +\langle w,\cA v\rangle$, so that $H$ is affine and skew-symmetric. Then for every $y\in Y$ it holds
\begin{equation}\label{coreq}
\epsilonsad(y|f,V,W) \leq\epsilonvi(y|H,Y).
\end{equation}
\end{corollary}
{\bf Proof.} Consider accuracy certificate $\cC^1=\{y_1=y,\lambda_1=1,H(y_1)\}$; for this certificate, $\widehat{y}$ as defined in Lemma \ref{lem2} is just $y$. Therefore, by Lemma \ref{lem2}, $\Res(\cC^1|Y)=\epsilonvi(y|H,Y)$. This equality combines with Lemma \ref{lem1} to imply (\ref{coreq}). \qed

\section{Representations of Monotone Operators}\label{secrepr}
\subsection{Outline} To explain the origin of the developments to follow, let us summarize the approach  to solving convex minimization problems on domains given by Linear Minimization Oracles (LMOs), developed in \cite{DualSubgr}. The principal ingredient of this approach is a {\sl Fenchel-type representation} of a convex function $f:X\to\bR$ defined on a convex subset $X$ of Euclidean space $E$; by definition, such a representation is
\begin{equation}\label{FTR}
f(x)=\min_{y\in Y}\left[\langle x,Ay+a\rangle -\psi(y)\right],
\end{equation}
where $Y$ is a convex subset of Euclidean space $F$ and $\psi:Y\to\bR$ is convex. Assuming for the sake of simplicity that $X$, $Y$ are compact and $\psi$ is continuously differentiable on $Y$, representation (\ref{FTR}) allows to associate with the {\sl primal problem}
$$
\Opt(P)=\min_{x\in X} f(x)\eqno{(P)}
$$
its dual
$$\Opt(D)=\max_{y\in Y}\left[f_*(y)=\min_{x\in X}\langle x,Ay+a\rangle -\psi(y)\right]\eqno{(D)}
$$
with the same optimal value. Observe that the first order information on the (concave) objective of $(D)$ is readily given by the first order information on $\psi$ and the information provided by an LMO for $X$. As a result, we can solve $(D)$ by, say, a proximal type First Order Method, provided that $Y$ is proximal-friendly. The crucial in this approach question of how to recover a good approximate solution to the problem of interest $(P)$ from the information collected when solving $(D)$ is addressed via the machinery of accuracy certificates \cite{NOR,DualSubgr}.
\par
In the sequel, we intend to apply a similar scheme to the situation where the role of $(P)$ is played by a variational inequality with monotone operator on a convex compact domain $X$ given by an LMO. Our immediate task is to outline {\sl informally}  what a Fenchel-type {\sl representation} of a monotone operator is and how we intend to use such a representation. To this end note that
$(P)$ and $(D)$ can be reduced to variational inequalities with monotone operators, specifically
\begin{itemize}
\item  the ``primal'' v.i. stemming from $(P)$. The domain of this v.i. is $X$,  and the operator is $f'(x)=Ay(x)+a$, where $y(x)$ is a maximizer of the function $\langle x,Ay\rangle -\psi(y)$ over $y\in Y$, or, which is the same, a (strong) solution to the v.i. given by the domain $Y$ and the monotone operator $y\mapsto G(y)-A^*x$, where $G(y)=\psi'(y)$;
\item the ``dual'' v.i. stemming from $(D)$.  The domain of this v.i. is $Y$, and the operator is $y\mapsto G(y)-A^*x(y)$, where $x(y)$ is  a minimizer of $\langle x,Ay+a\rangle$ over $x\in X$.
\end{itemize}
Observe that both operators in question are described in terms of a monotone operator $G$ on $Y$ and affine mapping $y\mapsto Ay+a:F\to E$; in the above construction $G$ was the gradient field of $\psi$, but the construction of the primal and the dual v.i.'s  makes sense whenever $G$ is  a monotone operator on $Y$ satisfying minimal regularity assumptions. The idea of the approach we are about to develop is as follows: in order to solve a v.i. with a monotone operator $\Phi$ and domain $X$ given by an LMO,
\begin{itemize}
\item[A.] We represent $\Phi$ in the form of $\Phi(x)=Ay(x)+a$, where $y(x)$ is a strong solution to the v.i. on $Y$ given by the operator $G(y)-A^*x$, $G$ being an appropriate monotone operator on $Y$.\\
    It can be shown that a desired representation always exists, but by itself existence does not help much -- we need the representation to be suitable for numerical treatment,  to be available in a ``closed computation-friendly form.''  We show that  ``computation-friendly'' representations of monotone operators admit a kind of {\sl fully algorithmic} calculus which, for all basic monotonicity-preserving operations, allows to get straightforwardly a desired representation of the result of an operation from the representations of the operands. In view of this calculus, ``closed analytic form'' representations, allowing to compute efficiently the values of monotone operators, automatically lead to required computation-friendly representations.
\item[B.] We use the representation from A to build the ``dual'' v.i. with domain $Y$ and the operator $\Theta(y)=G(y)-A^*x(y)$, with exactly the same $x(y)$ as above, that is, $x(y)\in\Argmin_{x\in X} \langle x,Ay+a\rangle$.  We shall see that $\Theta$ is monotone, and that usually there is a significant freedom in choosing $Y$; in particular, we typically can choose $Y$ to be proximal-friendly.
\item[C.] We solve the dual v.i. by an algorithm, like Mirror Descent or Mirror Prox, which produce necessary accuracy certificates.  We will see -- and this is our main result -- that such a certificate $\cC^N$  can be converted straightforwardly into a feasible solution $x^N$ to the v.i. of interest such that $\epsilonvi(x^N|\Phi,X)\leq\Res(\cC^N|Y)$. As a result, if the certificates in question are {\sl good}, meaning that the resolution of $\cC^N$ as a function of $N$ obeys the standard efficiency estimates of the algorithm used to solve the dual v.i., we solve the v.i. of interest with the same efficiency estimate as the one for the dual v.i. It remains to note that most of the existing first order algorithms for solving v.i.'s with monotone operators (various versions of polynomial time cutting plane algorithms, like the Ellipsoid method, Subgradient/Mirror Descent, and different bundle-level versions of Mirror Descent) indeed produce good accuracy certificates, see \cite{NOR,DualSubgr}.
\end{itemize}

\subsection{The construction}\label{construction}
\subsubsection{Situation}\label{secsit}
Consider the situation where we are given
\begin{itemize}
\item an affine mapping
$$
y\mapsto Ay+a: F\to E,
$$
where $E$, $F$ are Euclidean spaces;
\item a nonempty closed convex set $Y\subset F$;
\item a continuous monotone operator
$$
G(y):Y\to F
$$
which is {\sl good} w.r.t. $A,Y$, goodness meaning that the variational inequality $\VI(G(\cdot)-A^*x,Y)$ has a strong solution for every $x\in E$. Note that when $Y$ is convex compact, every continuous monotone operator on $Y$ is good, whatever be $A$;
\item a nonempty convex compact set $X$ in $E$.
\end{itemize}
These data give rise to two operators: ``primal'' $\Phi:X\to E$ which is monotone, and ``dual'' $\Psi:Y\to F$ which is antimonotone (that is, $-\Psi$ is monotone).
\subsubsection{Primal monotone operator}
The {\sl primal operator} $\Phi: E\to E$ is defined by
\begin{equation}\label{eq11}
\Phi(x)=Ay(x)+a: \,y(x)\in Y, \;\langle A^*x-G(y(x)),y(x)-y\rangle\geq 0,\,\,\forall y\in Y.
\end{equation}
Observe that {\sl required $y(x)$ do exist: these are just strong solutions to the variational inequality  given by the monotone operator $G(y)-A^*x$ and the domain $Y$.}\par
Now, with $x',x''\in E$,  setting $y(x')=y'$, $y(x'')=y''$, so that $y',y''\in Y$, we have
$$
\begin{array}{l}
\langle\Phi(x')-\Phi(x''),x'-x''\rangle =\langle Ay'-Ay'',x'-x''\rangle =
\langle y'-y'',A^*x'-A^*x''\rangle\\
 =\langle y'-y'',A^*x'\rangle + \langle y''-y',A^*x''\rangle\\
=\langle y'-y'',A^*x'-G(y')\rangle +\langle y''-y',A^*x''-G(y'')\rangle +\langle G(y'),y'-y''\rangle
+\langle G(y''),y''-y'\rangle \\
=\underbrace{\langle y'-y'',A^*x'-G(y')\rangle}_{\geq0\hbox{\scriptsize\rm\ due to $y'=y(x')$}}+\underbrace{\langle y''-y',A^*x''-G(y'')\rangle}_{\geq0 \hbox{\scriptsize\rm\ due to $y''=y(x'')$}}
+\underbrace{\langle G(y')-G(y''),y'-y''\rangle}_{\geq0\hbox{\scriptsize\rm\ since $G$ is monotone}}\geq0.
\end{array}
$$
Thus, {\sl $\Phi(x)$ is monotone.} We call (\ref{eq11}) {\sl a representation} of the monotone operator $\Phi$, and the data $F,A,a,y(\cdot),G(\cdot),Y$ -- the {\sl data} of the representation. We also say that these data {\sl represent} $\Phi$. \par
Given a convex domain $X\subset E$ and a monotone operator $\bar{\Phi}$ on this domain, we say that  data $F$, $A$, $a$, $y(\cdot)$, $G(\cdot)$, $Y$ of the above type {\sl represent $\bar{\Phi}$ on $X$}, if the monotone operator $\Phi$ represented by these data coincides with $\bar{\Phi}$ on $X$.
\subsubsection{Dual operator}
Operator $\Psi:Y\to F$ is given by
\begin{equation}\label{dualop}
\Psi(y)=A^*x(y)-G(y): \,x(y)\in X,\; \langle Ay+a,x(y)-x\rangle\leq 0,\,\,\forall x\in X
\end{equation}
(in words: $\Psi(y)=A^*x(y)-G(y)$, where $x(y)$ minimizes $\langle Ay+a,x\rangle$ over $x\in X$). This operator clearly is antimonotone, as the sum of two antimonotone operators $-G(y)$ and $A^*x(y)$; antimonotonicity of the latter operator stems from the fact that it is obtained from the antimonotone operator $\psi(z)$ -- a section of the superdifferential $\Argmin_{x\in X} \langle z,x\rangle$ of a concave function -- by affine substitution of variables: $A^*x(y)=A^*\psi(Ay+a)$, and this substitution preserves antimonotonicity.
\begin{remark} Note that computing the value of $\Psi$ at a point $y$ reduces to computing $G(y)$, $Ay+a$, a single call to the Linear Minimization Oracle for $X$ to get $x(y)$, and computing $A^*x(y)$.
\end{remark}
\subsection{Calculus of representations}
\subsubsection{Multiplication by nonnegative constants} Let $F,A,a,y(\cdot),G(\cdot),Y$ represent a monotone operator $\Phi:E\to E$:
$$
\Phi(x)=Ay(x)+a: y(x)\in Y\ \mbox{ and }\ \langle A^*x-G(y(x)),y(x)-y\rangle \geq 0\,\,\forall y\in Y.
$$
For $\lambda\geq0$, we clearly have
$$
\lambda\Phi(x)=[\lambda A]y(x)+[\lambda a]: \langle [\lambda A]^*x-[\lambda G(y(x))],y(x)-y\rangle\geq 0\forall y\in Y,
$$
that is, a representation of $\lambda\Phi$ is given by $F,\lambda A,\lambda a,y(\cdot),\lambda G(\cdot),Y$; note that the operator $\lambda G$ clearly is good w.r.t. $\lambda A,Y$, since $G$ is good w.r.t. $A,Y$.
\subsubsection{Summation} Let $F_i,A_i,a_i,y_i(\cdot),G_i(\cdot),Y_i$, $1\leq i\leq m$, represent monotone operators  $\Phi_i(x):E\to E$:
$$
\Phi_i(x)=A_iy_i(x)+a_i: \;y_i(x)\in Y_i\ \mbox{ and }\ \langle A_i^*x-G_i(y_i(x)),y_i(x)-y_i\rangle\geq0\,\,\forall y_i\in Y_i.
$$
Then
$$
\begin{array}{l}
\sum_i\Phi_i(x)=[A_1,...,A_m][y_1(x);...;y_m(x)]+[a_1;...;a_m],\\
y(x):=[y_1(x);...;y_m(x)]\in Y:=Y_1\times...\times Y_m,\\
\langle [A_1,...,A_m]^*x -[G_1(y_1(x));...;G_m(y_m(x))],[y_1(x);...;y_m(x)]-[y_1;...;y_m]\rangle\\
=\sum_i\langle A_i^*x-G_i(y_i(x)),y_i(x)-y_i\rangle \geq0\,\;\forall
y=[y_1;...;y_m]\in Y,
\end{array}
$$
so that the data
\bse
&&F=F_1\times...\times F_m,\;A=[A_1,...,A_m],\;[a_1;...;a_m],\\
&&y(x)=[y_1(x);...;y_m(x)],\;G(y)=[G_1(y_1);...;G_m(y_m)],\;Y=Y_1\times...\times Y_m
\ese represent $\sum_i\Phi_i(x)$. Note that the operator $G(\cdot)$ clearly is good since $G_1,...,G_m$ are so.
\subsubsection{Affine substitution of argument} Let $F,A,a,y(\cdot),G(\cdot),Y$ represent $\Phi:E\to E$, let $H$ be a Euclidean space and $h\mapsto Qh+q$ be an affine mapping from $H$ to $E$. We have
$$
\begin{array}{l}
\widehat{\Phi}(h):=Q^*\Phi(Qh+q)=Q^*(Ay(Qh+q)+a): \\
y(Qh+q)\in Y \mbox{ and }  \langle A^*[Qh+q]-G(y(Qh+q)),y(Qh+q)-y\rangle \geq0\,\,\forall y\in Y\\
\Rightarrow \hbox{with $\widehat{A}=Q^*A$, $\widehat{a}=Q^*a$, $\widehat{G}(y)=G(y)-A^*q$, $\widehat{y}(h)=y(Qh+q)$ we have}\\
\widehat{\Phi}(h)=\widehat{A}\widehat{y}(h)+\widehat{a}: \widehat{y}(h)\in Y\ \mbox{ and }\ \langle \widehat{A}^*h-\widehat{G}(\widehat{y}(h)),\widehat{y}(h)-y\rangle \geq0\,\;\forall y\in Y,\\
\end{array}
$$
that is, $F,\widehat{A},\widehat{a},\widehat{y}(\cdot),\widehat{G}(\cdot),Y$ represent $\widehat{\Phi}$. Note that $\widehat{G}$ clearly is good since $G$ is so.
\subsubsection{Direct sum} Let $F_i,A_i,a_i,y_i(\cdot),G_i(\cdot),Y_i$, $1\leq i\leq m$, represent monotone operators $\Phi_i(x_i):E_i\to E_i$. Then
$$
\begin{array}{l}
\Phi(x:=[x_1;...;x_m])=[\Phi_1(x_1);...;\Phi_m(x_m)]=\Diag\{A_1,...,A_m\}y(x)+[a_1;...;a_m]:\\
 y(x):=[y_1(x_1);...;y_m(x_m)]\in Y:=Y_1\times...\times Y_m
\ \mbox{ and }\ \\
\langle\Diag\{A^*_1,...,A^*_m\}[x_1;...;x_m] -[G_1(y_1(x_1));...;G_m(y_m(x_m))],y(x)-[y_1;...;y_m]\rangle\\
 =\sum_i\langle A_i^*x-G_i(y_i(x_i)),y_i(x_i)-y_i\rangle \geq0\,\;\forall y=[y_1;...;y_m]\in Y,\\
\end{array}
$$
so that
\bse
&F=F_1\times...\times F_m, \;A=\Diag\{A_1,...,A_m\},\;a=[a_1;...;a_m],\;y(x)=[y_1(x_1);...;y_m(x_m)],\\
&G(y=[y_1;...;y_m])=[G_1(y_1);...;G_m(y_m)],\; Y=Y_1\times...\times Y_m
 \ese
 represent $\Phi:\;E_1\times...\times E_m\to E_1\times...\times E_m$. Note that $G$ clearly is good since $G_1,...,G_m$ are so.
\subsubsection{Representing affine monotone operators}\label{sec335}
Consider an affine monotone operator on a Euclidean space $E$:
\begin{equation}\label{AffinePrimal}
\begin{array}{c}
\Phi(x)=Sx+a:E\to E\\
\left[S: \langle x,Sx\rangle\geq 0\,\forall x\in E\right]\\
\end{array}
\end{equation}
Its Fenchel-type representation on a convex compact set $X\subset E$ is readily given by the data $F=E$, $A=S$, $G(y)=S^*y:F\to F$ (this operator indeed is monotone), $y(x)=x$ and $Y$ being either the entire $F$, or (any) compact convex subset of $E=F$ which contains $X$; note that $G$ clearly is good w.r.t. $A,F$, same as is good w.r.t. $A,Y$ when $Y$ is compact. To check that the just defined $F,A,a,y(\cdot),G(\cdot),Y$ indeed represent
$\Phi$ on $X$, observe that when $x\in X$, $y(x)=x$ belongs to $Y\supset X$ and clearly satisfies the relation $0\leq \langle A^*x-G(y(x)),y(x)-y\rangle \geq0$ for all $y\in Y$ (see
(\ref{eq11})), since $A^*x-G(y(x))=S^*x-S^*x=0$. Besides this, for $x\in X$ we have
$$
Ay(x)+a=Sx+a=\Phi(x),
$$
as required for a representation. The dual antimonotone operator associated with this representation of $\Phi$ on $X$ is
\begin{equation}\label{AffineDual}
\Psi(y)=S^*[x(y)-y],\,\, x(y)\in\Argmin_{x\in X} \langle x,Sy+a\rangle.
\end{equation}
\subsubsection{Representing gradient fields} Let $f(x)$ be a convex function given by Fenchel-type representation
\begin{equation}\label{ftfr}
f(x)=\max\limits_{y\in Y}\left\{\langle x,Ay+a\rangle -\psi(y)\right\},
\end{equation}
where $Y$ is a convex compact set in Euclidean space $F$, and $\psi(\cdot):F\to \bR$ is a continuously differentiable convex function. Denoting by $y(x)$ a maximizer of $\langle x,Ay\rangle-\psi(y)$ over $y$, observe that
$$
\Phi(x):=Ay(x)+a
$$
is a subgradient field of $f$, and that this monotone operator is given by a representation with the data $F,A,a,y(\cdot),G(\cdot):=\nabla\psi(\cdot),Y$; $G$ is good since $Y$ is compact.
\section{Main result}\label{secmain} Consider the situation described in section \ref{construction}. Thus, we are given Euclidean space $E$, a convex compact set $X\subset E$ and a monotone operator $\Phi:E\to E$ represented according to (\ref{eq11}), the data being $F,A,a,y(\cdot),G(\cdot),Y$. We denote by $\Psi:Y\to F$ the dual (antimonotone) operator induced by the data $X,A,a,y(\cdot),G(\cdot)$, see (\ref{dualop}). Our goal is to solve variational inequality given by $\Phi$, $X$, and our main observation is that {\sl a good accuracy certificate for the variational inequality given by $(-\Psi,Y)$ induces an equally good solution to the variational inequality given by $(\Phi,X)$.} The exact statement is as follows.
\begin{theorem}\label{propmain} Let $X\subset E$ be a convex compact set and $\Phi: X\to E$ be a monotone operator represented on $X$, in the above sense,  by data $F,A,a,y(\cdot),G(\cdot),Y$. Let also $\Psi:Y\to F$ be the antimonotone operator as defined above by the data $X,A,a,y(\cdot),G(\cdot)$. Let, finally, $$\cC^N=\{y_t,\lambda_t,-\Psi(y_t)\}_{t=1}^N$$ be an accuracy certificate associated with the  monotone operator $[-\Psi]$ and $Y$. Setting
$$
x_t=x(y_t)
$$
(these points are byproducts of computing $\Psi(y_t)$, $1\leq t\leq N$) and
$$
\widehat{x}=\sum_{t=1}^N\lambda_t x_t\;(\in X),
$$
we ensure that
\begin{equation}\label{ensure}
\epsilonvi(\widehat{x}|\Phi,X)\leq \Res\left(\cC^N|Y(X)\right),\,\,Y(X):=\{y(x):x\in X\}\subset Y.
\end{equation}
When $\Phi(x)=a+Sx$, $x\in X$, with skew-symmetric $S$, we have also
\begin{equation}\label{ensureskew}
\Res(\{x_t,\lambda_t,\Phi(x_t)\}_{t=1}^N|X)\leq \Res\left(\cC^N|Y(X)\right).
\end{equation}
\end{theorem}
\par
In view of Theorem \ref{propmain}, given a representation of the monotone operator $\Phi$ participating in the v.i. of interest $\VI(\Phi,X)$, we can reduce solving the v.i. to solving the {\sl dual v.i.} $\VI(-\Psi,Y)$ by an algorithm producing good accuracy certificates.
Below we discuss in details the situation when the latter algorithm is either Mirror Descent (MD) \cite[Chapter 5]{chapters}, or Mirror Prox (MP)  (\cite{MP}, see also \cite[Chapter 6]{chapters} and \cite{Acta}).
\par
Theorem \ref{propmain} may be extended to the situation where the relationships \rf{dualop}, defining the dual operator $\Psi$ hold only approximately. We present here the following slight extension of the main result:
\begin{theorem}\label{propmainnew} Let $X\subset E$ be a convex compact set and $\Phi: X\to E$ be a monotone operator represented on $X$, in the sense of section \ref{secrepr},  by data $F,A,a,y(\cdot),G(\cdot),Y$. Given a positive integer $N$, sequences $y_t\in Y$, $x_t\in X$, $1\leq t\leq N$, and nonnegative reals $\lambda_t$, $1\leq t\leq N$, summing up to 1, let us set
\begin{equation}\label{epsilonetc}
\epsilon=\Res(\{y_t,\lambda_t,G(y_t)-A^*x_t\}_{t=1}^N|X)=\sup_{z\in Y(X)}\sum_{t=1}^N\lambda_t\langle G(y_t)-A^*x_t,y_t-z\rangle,
\end{equation}
and
\[
\widehat{x}=\sum_{t=1}^N \lambda_t x_t\;(\in X).
\]
Then
\begin{equation}\label{bottomline}
\epsilonvi(\widehat{x}|\Phi,X) \leq\epsilon+\sup_{x\in X}\sum_{t=1}^N\lambda_t\langle Ay_t+a,x_t-x\rangle.
\end{equation}
\end{theorem}
Proofs of Theorems \ref{propmain} and \ref{propmainnew} are given in section \ref{sec:ensure}.

\subsection{Mirror Descent and Mirror Prox algorithms}
\paragraph{Preliminaries.} Saddle Point MD and MP are  algorithms for solving convex-concave saddle point problems and variational inequalities with monotone operators\footnote{MD algorithm originates from \cite{NYu,NemVI}; its modern proximal form was developed in \cite{BeckTeb}. MP was proposed in \cite{MP}. For the most present exposition of the algorithms, see \cite[Chapters 5,6]{chapters} and \cite{Acta}.}. The algorithms are of {\sl proximal type}, meaning that in order to apply the algorithm to a v.i. $\VI(H,Y)$, where $Y$ is a nonempty closed convex set in Euclidean space $E_y$ and $H$ is a monotone operator on $Y$, one needs to equip $E_y$ with a norm $\|\cdot\|$, and $Y$ - with a continuously differentiable {\sl distance generating function} (d.-g.f.) $\omega(\cdot):Y\to\bR$ {\sl compatible} with $\|\cdot\|$, meaning that  $\omega$ is strongly convex, modulus 1, w.r.t. $\|\cdot\|$. We call  $\|\cdot\|,\omega(\cdot)$ {\sl proximal setup} for $Y$. This setup gives rise to
\begin{itemize}
\item {\sl $\omega$-center} $y_\omega=\argmin_{y\in Y}\omega(y)$ of $Y$,
\item {\sl Bregman distance}
$$
V_y(z)=\omega(z)-\omega(y)-\langle\omega'(y),z-y\rangle \geq {\half}\|z-y\|^2
$$
where the concluding inequality is due to strong convexity of $\omega$,
\item {\sl $\omega$-size of a nonempty subset $Y'\subset Y$}
$$
\Omega[Y']=\sqrt{2\left[\max_{y'\in Y'}\omega(y') -\min_{y\in Y}\omega(y)\right]}.
$$
Due to the origin of $y_\omega$, we have $V_{y_\omega}(y)\leq \half\Omega^2[Y']$ for all $y\in Y'$, implying that $\|y-y_\omega\|\leq \Omega$ for all $y\in Y'$.
\end{itemize}
Given $y\in Y$, the {\sl prox-mapping} with center $u$ is defined as
$$
\Prox_y(\zeta)=\argmin_{z\in Y}\left[V_y(z)+\langle \zeta,z\rangle\right]=\argmin_{z\in Y}\left[\omega(z)+\langle \zeta-\omega'(y),z\rangle\right]:\,E_y\to Y.
$$
\paragraph{The algorithms.} Let $Y$ be a nonempty closed convex set in Euclidean space $E_y$, $H=\{H_t:Y\to E_y\}_{t=1}^\infty$ be a sequence of vector fields, and $\|\cdot\|,\omega(\cdot)$ be a proximal setup for $Y$.
As applied to $(H,Y)$, MD  is the recurrence
\begin{equation}\label{MDalg}
\begin{array}{rcl}
y_1&=&y_\omega;\\
y_t&\mapsto& y_{t+1}=\Prox_{y_t}(\gamma_t H_t(y_t)),\,t=1,2,...
\end{array}
\end{equation}
MP is the recurrence
\begin{equation}\label{MPalg}
\begin{array}{rcl}
y_1&=&y_\omega;\\
y_t&\mapsto& z_t=\Prox_{y_t}(\gamma_t H_t(y_t)) \mapsto y_{t+1}=\Prox_{y_t}(\gamma_t H_t(z_t)),\,t=1,2,...
\end{array}
\end{equation}
In both MD and MP, $\gamma_t>0$ are stepsizes. The most important to us properties of these recurrences are as follows.
\begin{proposition}\label{propMD}
For $N=1,2,...$, consider the accuracy certificate
\[\cC^N=\bigg\{y_t\in Y,\;\lambda^N_t:={\gamma_t}\left[\sum_{\tau=1}^N\gamma_\tau\right]^{-1},\;H_t(y_t)\bigg\}_{t=1}^N,
\] associated with {\rm (\ref{MDalg})}.
Then for every  $Y'\subset Y$ one has
\begin{equation}\label{thenMD}
\Res(\cC^N|Y') \leq {\Omega^2[Y']+\sum_{t=1}^N \gamma_t^2\|H_t(y_t)\|_*^2\over 2\sum_{t=1}^N\gamma_t},
\end{equation}
where $\|\cdot\|_*$ is the norm conjugate to $\|\cdot\|$: $
\|\xi\|_*=\max_{\|x\|\leq1} \langle \xi,x\rangle.
$
\par\noindent
In particular, if
\begin{equation}\label{ifbounded}
\forall (y\in Y,t): \|H_t(y)\|_*\leq M
\end{equation}
with some finite $M\geq0$, then, given $Y'\subset Y$, $N$ and setting
\begin{equation}\label{MDgammasatisfy}
(a): \gamma_t={\Omega[Y']\over M\sqrt{N}},\,1\leq t\leq N,\hbox{\ or\ }(b): \gamma_t={\Omega[Y']\over \|H_t(y_t)\|_*\sqrt{N}},\,1\leq t\leq N,
\end{equation}
one has
\begin{equation}\label{MDrate}
\Res(\cC^N|Y')\leq {\Omega[Y'] M\over\sqrt{N}}.
\end{equation}
\end{proposition}
\begin{proposition}\label{propMP}
For $N=1,2,...$, consider the accuracy certificate
\[\cC^N=\bigg\{z_t\in Y,\;\lambda^N_t:={\gamma_t}\left[\sum_{\tau=1}^N\gamma_\tau\right]^{-1},\;H_t(z_t)\bigg\}_{t=1}^N,
\] associated with {\rm (\ref{MPalg})}.
Then, setting
\begin{equation}\label{dt}
d_t=\gamma_t\langle H_t(z_t),z_t-y_{t+1}\rangle -V_{y_t}(y_{t+1}),
\end{equation}
we have for every $Y'\subset Y$
\be\label{thenMP1}
\Res(\cC^N|Y') &\leq& {\half\Omega^2[Y']+\sum_{t=1}^N d_t\over\sum_{t=1}^N\gamma_t}\\
d_t&\leq&{\half}\left[\gamma_t^2\|H_t(z_t)-H_t(y_t)\|_*^2-\|y_t-z_t\|^2\right],
\ee{thenMP2}
where $\|\cdot\|_*$ is the norm conjugate to $\|\cdot\|$.
In particular, if
\begin{equation}\label{ifsmooth}
\forall (y,y'\in Y,t): \|H_t(y)-H_t(y')\|_*\leq L\|y-y'\|+M
\end{equation}
with some finite $L\geq0$, $M\geq0$, then given $Y'\subset Y$, $N$ and setting
\begin{equation}\label{gammasatisfy}
\gamma_t={1\over\sqrt{2}}\min\left[{1\over L},{\Omega[Y']\over M\sqrt{N}}\right],\,1\leq t\leq N,
\end{equation}
one has
\begin{equation}\label{MPrate}
\Res(\cC^N|Y')\leq {1\over\sqrt{2}}\max\left[{\Omega^2[Y']L\over N},{\Omega[Y'] M\over\sqrt{N}}\right].
\end{equation}
\end{proposition}
To make the text self-contained, we provide the proofs of these known results in the appendix.
\subsection{Intermediate summary}
Theorem \ref{propmain} combines with Proposition \ref{propMD} to imply the following claim:
\begin{corollary}\label{corbas}
In the situation of Theorem \ref{propmain}, let $y_1,...,y_N$ be the trajectory of $N$-step MD as applied to the stationary sequence $H=\{H_t(\cdot)=-\Psi(\cdot)\}_{t=1}^\infty$ of vector fields, and let $x_t=x(y_t)$, $t=1,...,N$. Then, setting $\lambda_t={\gamma_t\over\sum_{\tau=1}^N\gamma_\tau}$, $1\leq t\leq N$, we ensure that
\begin{equation}\label{then17}
\epsilonvi(\underbrace{\sum_{t=1}^N\lambda_tx_t}_{\widehat{x}}|\Phi,X)\leq \Res(\underbrace{\{y_t,\lambda_t,-\Psi(y_t)\}_{t=1}^N}_{\cC^N}|Y(X))\leq
{\Omega^2[Y(X)]+\sum_{t=1}^N \gamma_t^2\|\Psi(y_t)\|_*^2\over 2\sum_{t=1}^N\gamma_t}.
\end{equation}
In particular, assuming
$$
M=\sup_{y\in Y}\|\Psi(y)\|_*
$$
finite and specifying $\gamma_t$, $1\leq t\leq N$, according to {\rm (\ref{MDgammasatisfy})} with $Y'=Y(X)$, we ensure that
\begin{equation}\label{sqrtNefficiency}
\epsilonvi(\widehat{x}|\Phi,X)\leq \Res(\cC^N|Y(X))\leq {\Omega[Y(X)] M\over\sqrt{N}}.
\end{equation}
When $\Phi(x)=Sx+a$, $x\in X$, with a skew-symmetric $S$, $\epsilonvi(\widehat{x}|\Phi,X)$ in the latter relation can be replaced with
$\Res(\{x_t,\lambda_t,\Phi(x_t)\}_{t=1}^N|X)$.
\end{corollary}
In the sequel, we shall refer to the implementation of our approach presented in Corollary \ref{corbas} as to our {\sl basic} scheme.
\section{Modifications in Affine case}\label{secmodif}
In this section, we present some modifications of the proposed approach  as applied to the case of v.i. $\VI(\Phi,X)$ with {\sl affine} monotone operator $\Phi$ and LMO-represented convex compact domain $X$. While the worst-case complexity bounds for the modified scheme are similar to the ones stated in Corollary \ref{corbas}, there are reasons to believe that in practice the modified scheme could outperform the basic one.

\subsection{Situation}
In the rest of this section, we consider the case when the monotone operator $\Phi$ is affine:
$$
\Phi(x)=Sx+a:E\to E\eqno{[S: \langle Sx,x\rangle \geq 0\,\forall x\in E]}
$$
and our goal is to solve $\VI(\Phi,X)$, where $X$ is a convex compact subset of $E$; w.l.o.g. we assume that $0\in X$. We suppose that $\Phi$  is given by an {\sl affine} Fenchel-type representation, that is, a representation with data \begin{equation}\label{thenewdata}
F,\;A,\;a,\;G(y):=Gy,\;y(x):=Bx,\;Y=F,
\end{equation}
where
\begin{enumerate}
\item $F$ is a Euclidean space, $y\mapsto Ay+a$ is an affine mapping from $F$  to $E$;
\item $y\mapsto Gy:F\to F$ is a linear monotone mapping, and $x\to Bx: E\to F$ is a linear mapping such that
\begin{equation}\label{suchthat23}
\begin{array}{lrcllrcl}
(a)&GB&=&A^*,&
(b)&AB&=&S;\\
\end{array}
\end{equation}
\item $Y=F$.
\end{enumerate}
Note that (\ref{suchthat23}.$a$) implies that setting $y(x)=Bx$, $y(x)$ is a strong solution to the v.i. associated with the operator $G(y)-A^*x$ and $Y=F$, while (\ref{suchthat23}.$b$) says that $Ay(x)+a=\Phi(x)$, $x\in E$, as required in the definition \rf{eq11} of a Fenchel-type representation of a monotone operator.
\subsection{Strategy}\label{sec:strategy}
\subsubsection{Preliminaries} We intend to get an approximate solution to $\VI(\Phi,X)$ by applying MP to a properly built sequence $H=\{H_t(\cdot)\}$ of vector fields on $F$. Let us fix a  proximal setup $\|\cdot\|,\omega(\cdot)$ for $Y=F$; w.l.o.g., we assume that the $\omega$-center $\argmin_F\omega(\cdot)$ of $F$ is the origin, that is, $\omega'(0)=0$. Let $L$ be the operator norm of the mapping $y\mapsto G(y):=Gy:F\to F$ from $\|\cdot\|$ to $\|\cdot\|_*$, so that
$$
\forall y\in F: \|Gy\|_*\leq L\|y\|,
$$
or, equivalently, $\langle z,Gy\rangle \leq \|z\|\|y\|$ for all $z,y\in F$.
In the sequel, we set
$$
\gamma={L}^{-1}.
$$
\subsubsection{The construction}\label{sec:constr5}
We intend to build $H_t(\cdot)$ recursively, according to the recurrence
\begin{equation}\label{recurrence}
\begin{array}{rcl}
y_1&=&0\\
y_t&\mapsto& x_t\in X\mapsto H_t(v) = Gv-A^*x_t [\equiv G(v)-A^*x_t]\mapsto\\
z_t&=&\Prox_{y_t}(\gamma H_t(y_t))\mapsto y_{t+1} =\Prox_{y_t}(\gamma H_t(z_t)).\\
\end{array}
\end{equation}
Note that independently of the choice of $x_t\in E$, we have
\[
\|H_t(v)-H_t(v')\|_*\leq L\|v-v'\|.
\]
Now the relationships of the MP recurrence imply that (see (\ref{verylast-1}) and (\ref{verylast}))
 \begin{equation}\label{whence101}
\forall z\in F: \gamma\langle H_t(z_t), z_t-z\rangle \leq V_{y_t}(z)-V_{y_{t+1}}(z).
\end{equation}
The essence of the matter is how we update the vectors $x_t$; this is the issue we consider next.
\paragraph{Functions $f_y(\cdot)$.} Given $y\in F$, let us set
\begin{equation}\label{fyofx}
f_y(x)=
\gamma\langle a,x\rangle+\max_{z\in F} \left[\langle z,\gamma[A^*x-Gy]\rangle -V_y(z)\right]
\end{equation}
Since $\omega(\cdot)$ is strongly convex on $F$, the function $f_y(\cdot)$ is well defined on $E$; $f_y$ is convex as the supremum of a family of affine functions of $x$. Moreover, it is well known that in fact $f_y(\cdot)$ possesses Lipschitz continuous gradient.
Specifically, let $\|\cdot\|_E$ be a norm on $E$, $\|\cdot\|_{E,*}$ be the norm conjugate to $\|\cdot\|_E$, and let $\cL$ be the norm of the linear mapping $y\to Ay:F\to E$ from the norm $\|\cdot\|$ on $F$ to the norm $\|\cdot\|_{E,*}$ on $E$, so that
\[
\langle Ay,x\rangle \leq \cL\|y\|\|x\|_{E}\;\;\forall (y\in F,x\in E),
\]
or, what is the same,
\begin{equation}\label{calLnorm}
\begin{array}{rcl}
\|Ay\|_{E,*}&\leq& \cL\|y\|\,\,\forall y\in F\\
\|A^*x\|_*&\leq& \cL \|x\|_E\,\,\forall x\in E.
\end{array}
\end{equation}
\begin{lemma}\label{lemmaholdstrue} Let
$
z_y(\zeta)=\Prox_y(\zeta):F\to Y.$
Function  $f_y(\cdot)$ is continuously differentiable with the gradient
\begin{equation}\label{eqgradf}
\nabla f_y(x)=\gamma Az_y(\gamma[Gy-A^*x])+\gamma a,\\
\end{equation}
and this gradient is Lipschitz continuous:
\begin{equation}\label{LipContGrad}
\|\nabla f(x')-\nabla f(x'')\|_{E,*}\leq (\gamma\cL)^2\|x'-x''\|_E\,\,\forall x',x''\in E.
\end{equation}
\end{lemma}
For proof, see section \ref{Proofoflemmaholdstrue}.
\paragraph{Updating $x_t$'s, preliminaries.} Observe, first, that
when summing up inequalities (\ref{whence101}), we get
\begin{equation}\label{weget17}
\Res(\{y_t,\lambda_t=N^{-1},H_t(z_t)\}_{t=1}^N|Y(X))\leq {1\over 2\gamma N}\Omega^2[Y(X)]={\Omega^2[Y(X)]L\over 2N},\quad Y(X)=BX.
\end{equation}
Second, for any $x_t\in X$, $1\leq t\leq N$, we have $\widehat{x}={1\over N}\sum_{t=1}^Nx_t\in X$. Further, invoking (\ref{bottomline}) with $\lambda_t=N^{-1}$, $1\leq t\leq N$, and $z_t$ in the role of $y_t$ (which by (\ref{weget17}) allows to set $\epsilon={\Omega^2[Y(X)]L\over 2N}$), we get
\begin{equation}\label{weget18}
\begin{array}{rcl}
\epsilonvi(\widehat{x}|\Phi,X)&=&\max\limits_{{x}\in X}\langle \Phi({x}),\widehat{x}-{x}\rangle\\
&\leq& {L\Omega^2[Y(X)]\over 2N}+
\max_{{x}\in X}{1\over N}\sum_{t=1}^N\langle Az_t+a,x_t-{x}\rangle\\
&=& {L\Omega^2[Y(X)]\over 2N}+\max\limits_{{x}\in X}{L\over N}\sum_{t=1}^N\langle \nabla f_{y_t}(x_t),x_t-{x}\rangle\\
\end{array}
\end{equation}
(we have used (\ref{eqgradf}) and have taken into account that $z_t=z_{y_t}(\gamma[Gy_t-A^*x_t])$, see (\ref{eqgradf}) and (\ref{recurrence}); recall that $\gamma=1/L$). Note that so far our conclusions were  independent on how $x_t\in X$ are selected.
\par
Relation (\ref{weget18}) implies that {\sl when $x_t$ is a minimizer of $f_{y_t}(\cdot)$ on $X$, we have $\langle \nabla f_{y_t}(x_t),x_t-{x}\rangle \leq0$ for all ${x}\in X$, and with this ``ideal'' for our purposes choice of $x_t$,
{\rm(\ref{weget18})} would imply
$$
\epsilonvi(\widehat{x}|\Phi,X)\leq {L\Omega^2[Y(X)]\over 2N},
$$
which is an $O(1/N)$ efficiency estimate}, much better that the $O(1/\sqrt{N})$-efficiency estimate {\rm (\ref{sqrtNefficiency})}.
\paragraph{Updating $x_t$'s, CGA implementation.}
Of course, we cannot simply specify $x_t$ as a point from $\Argmin_X f_{y_t}(x)$, since this would require solving precisely at every step of the MP recurrence (\ref{recurrence}) a large-scale convex optimization problem. What we indeed intend to do, is to solve this problem {\sl approximately.} Specifically, given $y_t$ (so that $f_{y_t}(\cdot)$ is identified),  we can apply the classical Conditional Gradient Algorithm (CGA) (which, as was explained in the introduction, is, basically, the only traditional algorithm capable to minimize a smooth convex function over an LMO-represented convex compact set) in order to generate an approximate solution $x_t$ to the problem $\min_{X} f_{y_t}(x)$ satisfying, for some prescribed $\epsilon>0$, the relation
\begin{equation}\label{deltat}
\delta_t:=\max_{x\in X} \langle \nabla f_{y_t}(x_t),x_t-x\rangle \leq \epsilon.
\end{equation}
By (\ref{weget18}), this course of actions implies the efficiency estimate
\begin{equation}\label{implyeffest}
\epsilonvi(\widehat{x}|\Phi,X)\leq {L\Omega^2[Y(X)]\over 2N}+L\epsilon.
\end{equation}
\subsubsection{Complexity analysis}\label{sec:complexanal}
 Let us equip $E$  with a norm $\|\cdot\|_E$, the conjugate norm being $\|\cdot\|_{E,*}$, and let $\cL$ be the operator norm of the mapping $y\mapsto Ay$ as defined in Lemma \ref{lemmaholdstrue}. Let, further, $R=R_E(X)$ be the radius of the smallest $\|\cdot\|_E$-ball, centered at the origin, which contains $X$.  Taking into account (\ref{LipContGrad})
and applying the standard results on CGA (see section \ref{sec:CGA}), for every $\epsilon\in(0,\cL R^2)$ it takes at most $O(1)\cL R_E^2(X)/\epsilon$ CGA steps to generate a point $x_t$ with $\delta_t\leq\epsilon$; here and below $O(1)$'s are absolute constants. Specifying $\epsilon$ as ${\Omega^2[Y(X)]\over 2N}$, (\ref{implyeffest}) becomes
$$
\epsilonvi(\widehat{x}|\Phi,X)\leq {L\Omega^2[Y(X)]\over N},
$$
while the computational effort to generate $\widehat{x}$ is dominated by the necessity to generate $x_1,...,x_N$, which amounts to the total of
 $$
 \cN(N)=O(1){\cL R_E^2(X)\over \Omega^2[Y(X)]}N^2
 $$
 CGA steps. The effort per step is dominated by the necessity to compute the vector $g=\nabla f_y(x)$, given $y\in F$, $x\in E$, and to minimize the linear form $\langle g,u\rangle$ over $u\in X$. In particular, to ensure $\epsilonvi(\widehat{x}|\Phi,X)\leq\epsilon$, the total number of CGA steps should be proportional to $1/\epsilon^2$. We see that in terms of the theoretical upper bound on the number of calls to the LMO for $X$ needed to get an $\epsilon$-solution, our current scheme has no advantages as compared to the MD-based approach analyzed in Corollary \ref{corbas}. We, however, may hope that in practice the outlined MP-based scheme can be better than our basic MD-based one, provided that we apply CGA in a ``smart'' way, e.g., use CGA with memory, see \cite{HJN}.
\section{Illustration}\label{secill}
\subsection{The problem.} We apply our construction to the following problem (``matrix completion with spectral norm fit''):\footnote{A more interesting for applications problem (cf. \cite{CP2010,CP2011,JKN2012}) would be
\[
\Opt=\min_{v\in\bR^{p_v\times q_v}}\left\{\|v\|_\nuc: \|\cA v - b\|_{2,2}\leq\delta\right\};
\]
applying the approach from \cite{BSPP}, this problem can be reduced to a ``small series'' of problems (\ref{fineq1}).}
\begin{equation}\label{fineq1}
\Opt(P)=\min_{v\in\bR^{p_v\times q_v}:\|v\|_\nuc\leq 1}\left[\overline{f}(v):=\|\cA v - b\|_{2,2}\right]
\end{equation}
where $\bR^{p\times q}$ is the space of $p\times q$ real matrices, $\|x\|_\nuc=\sum_i\sigma_i(x)$ is the nuclear norm on this space (sum of the singular values $\sigma_i(x)$ of $x$), $\|x\|_{2,2}=\max_i\sigma_i(x)$ is the spectral norm of $x$ (which is exactly the conjugate of the nuclear norm), and $\cA$ is a linear mapping from $\bR^{p_v\times q_v}$ to $\bR^{p_b\times q_b}$.
We are interested in the ``large-scale'' case, where the sizes of $p_v,q_v$ of $v$ are large enough to make the full singular value decomposition of a $p_v\times q_v$ matrix prohibitively time consuming, what seemingly rules out
the possibility to solve (\ref{fineq1}) by proximal type First Order algorithms. We assume, at the same time, that
computing the leading singular vectors and the leading singular value of a $p_v\times q_v$ or a $p_b\times q_b$ matrix (which, computationally, is by far easier task than finding full singular value decomposition) still can be carried out in reasonable time.
\subsubsection{Processing the problem}
We rewrite (\ref{fineq1}) as a bilinear saddle point problem
\begin{equation}\label{saddlepoint}
\begin{array}{c}
\Opt(P)=\min\limits_{v\in V}\max\limits_{w\in W}\underbrace{\langle w,[\cA v-b]\rangle_\Fro}_{f(v,w)}
\\
V=\{v\in\bR^{p_v\times q_v}:\|v\|_\nuc\leq1\},\,\,W=\{w\in\bR^{p_b\times q_b}:\|w\|_\nuc\leq1\}\\
\end{array}
\end{equation}
(from now on $\langle\cdot,\cdot\rangle_\Fro$ stands for Frobenius inner product, and $\|\cdot\|_\Fro$ -- for the Frobenius norm on the space(s) of matrices).
The domain $X$ of the problem is the direct product of two unit nuclear norm balls; minimizing a linear form over this domain reduces to minimizing, given $\xi$ and $\eta$,  the linear forms $\Tr(v\xi^T)$, $\Tr(w\eta^T)$ over $\{v\in\bR^{p_v\times q_v}:\|v\|\leq1\}$, resp.,  $\{w\in\bR^{p_b\times q_b}:\|w\|\leq1\}$, which, in turn, reduces to computing the leading singular vectors and singular values of $\xi$ and $\eta$.
\par
The monotone operator associated with (\ref{saddlepoint}) is affine and skew-symmetric:
\[
\Phi(v,w)=\left[\nabla_vf(v,w);-\nabla_wf(v,w)\right]=\left[\cA^*w;-\cA v\right]+\left[0;b\right]:\;\underbrace{\bR^{p_v\times q_v}\times\bR^{p_b\times q_b}}_{E}\to E.
\]
 From now on we assume that {\sl $\cA$ is of spectral norm at most 1},  i.e.,
 \[
 \|\cA v\|_\Fro\leq \|v\|_\Fro,\;\; \forall v
 \]
  (this always can be  achieved by scaling).
\paragraph{Representing $\Phi$.}
We can represent the restriction of $\Phi$ on $X$ by the data
\begin{equation}\label{thedata}
\begin{array}{rcl}
F&=&\bR^{p_v\times q_v}\times \bR^{p_v\times q_v}\\
A y+a&=&[\xi;\cA \eta +b],\,y=[\xi;\eta]\in F \;(\xi\in \bR^{p_v\times q_v}, \eta\in \bR^{p_v\times q_v}), \\
G(\underbrace{[\xi;\eta]}_{y})&=&[-\eta;\xi]:\; F\to F\\
Y&=&\{y=[\xi;\eta]\in F: \;\|\xi\|_\Fro\leq1,\;\|\eta\|_\Fro\leq1\}\\
\end{array}
\end{equation}
Indeed, in the notation from section \ref{construction}, for $x=[v;w]\in X=\{[v;w]\in\bR^{p_v\times q_v}\times\bR^{p_b\times q_b}:\|v\|_\nuc\leq1,\|w\|_\nuc\leq1\}$, the solution $y(x)=[\xi(x);\eta(x)]$ to the linear system $A^*x=G(y)$ is given by
$\eta(x)=-v$, $\xi(x)=\cA^*w$, so that both components of $y(x)$ are of Frobenius norm $\leq1$ (recall that spectral norm of $\cA$ is $\leq1$), and therefore $y(x)\in Y$. Besides this,
$$
Ay(x=[v;w])+a=[\xi(x);\cA \eta(x)+b]=[\cA^*w;b-\cA v]=\Phi(v,w).
$$
We conclude that when $x=[v;w]\in X$, the just defined $y(x)$ meets all requirements from (\ref{eq11}), and thus the data $F,A,a,y(\cdot),G(\cdot),Y$ given by (\ref{thedata}) indeed represent the monotone operator $\Phi$ on $X$.
\paragraph{The dual operator}  $\Psi$ given by the data $F,A,a,y(\cdot),G(\cdot),Y$  is
\begin{equation}\label{dualpsi}
\begin{array}{l}
\Psi(\overbrace{[\xi;\eta]}^{y})=A^*x(y)-G(y)=[v(y)+\eta;\cA^*w(y)-\xi],\\
v(y)\in\Argmin\limits_{\|v\|_\nuc\leq 1}\langle v,\xi\rangle,\;\;w(y)\in\Argmin\limits_{\|w\|\leq1}\langle w,\cA \eta+b\rangle.
\end{array}
\end{equation}
\paragraph{Proximal setup.} We use the Euclidean proximal setup for $Y$, i.e., we equip the space $F$ embedding $Y$ with the Frobenius norm and take, as the d.-g.f. for $Y$, the function
$$
\omega(\xi,\eta) =\half\left[\|\xi\|_\Fro^2+\|\eta\|_\Fro^2\right]: F:=\bR^{p_v\times q_v}\times \bR^{p_v\times q_v}\to\bR,
$$
resulting in
$
\Omega[Y]=\sqrt{2}.
$
Furthermore, from (\ref{dualpsi}) and the fact that the spectral norm of $\cA$ is bounded by 1 it follows that
the monotone operator $\Theta(y)=-\Psi(y):Y\to F$ satisfies (\ref{ifbounded}) with $M=2\sqrt{2}$ and (\ref{ifsmooth}) with $L=0$ and $M=4\sqrt{2}$.
\paragraph{Remark.} Theorem \ref{propmain} combines with Corollary \ref{newcor} to imply that when converting an accuracy certificate $\cC^N$ for the dual v.i. $\VI(-\Psi,Y)$ into a feasible solution $\widehat{x}^N$ to the primal v.i. $\VI(\Phi,X)$, we ensure that
\begin{equation}\label{spacc}
\epsilonsad(\widehat{x}^N|f,V,W)\leq\Res(\cC^N|Y(X))\leq \Res(\cC^N|Y),
\end{equation}
with $f,V,W$ given by (\ref{saddlepoint}). In other words, in the representation $\widehat{x}^N=[\widehat{v}^N;\widehat{w}^N]$,
$\widehat{v}^N$  is a feasible solution to problem (\ref{fineq1}) (which is the primal problem associated with (\ref{saddlepoint})),
and $\widehat{w}^N$ is a feasible solution to the problem
\[
\Opt(D)=\max\limits_{w\in W} \min_{v\in  V}\langle w,\cA v -b\rangle= \max\limits_{w\in W}\left\{\underline{f}(w):=-\|A^*w\|_{2,2}-\langle b,w\rangle\right\},
\]
(which is the dual problem associated with (\ref{saddlepoint})) with the sum of non-optimalities, in terms of respective objectives, $\leq \Res(\cC^N|Y)$.
Computing $\underline{f}(\widehat{w})$ (which, together with computing $\overline{f}(\widehat{v})$, takes a single  call to LMO for $X$), we get a lower bound on
$\Opt(P)=\Opt(D)$ which certifies that $\overline{f}(\widehat{v}) -\Opt(P)\leq\Res(\cC^N|Y)$.
\subsection{Numerical illustration}
  Here we report on some numerical experiments with problem (\ref{fineq1}). In these experiments, we used $p_b=q_b=:m$, $p_v=q_v=:n$, with $n=2m$, and  the mapping $\cA$ given by
\begin{equation}\label{ourcA}
\cA v=\sum_{i=1}^k \ell_i v r_i^T,
\end{equation}
with generated at random  $m\times n$ factors $\ell_i,r_i$ scaled to get $\|\cA\|_{*}\approx1$. In all our experiments, we used $k=2$. Matrix $b$ in (\ref{fineq1}) was built as follows: we  generated at random $n\times n$ matrix $\bar{v}$ with $\|\bar{v}\|_\nuc$ less than (and close to) 1 and $\rank(\bar{v})\approx\sqrt{n}$, and took $b=\cA\bar{v}+\delta$, with randomly generated $m\times m$ matrix $\delta$ of spectral norm about 0.01.
\subsubsection{Experiments with the MD-based scheme}\label{sec:exp-MD}
\paragraph{Implementing the MD-based scheme.} In the first series of experiments, the dual v.i. $\VI(-\Psi,Y)$ is solved  by the MD algorithm with $N=512$ steps for all but the largest instance, where $N=257$ is used. The MD is applied to the stationary sequence $H_t\equiv -\Psi$, $t=1,2,...$, of vector fields. The stepsizes $\gamma_t$ are proportional, with coefficient of order of 1, to those given by (\ref{MDgammasatisfy}.$b$) with $\|\cdot\|\equiv\|\cdot\|_*=\|\cdot\|_\Fro$ and $\Omega[Y]=\sqrt{2}$ \footnote{As we have already mentioned, with our proximal setup, the $\omega$-size of $Y$ is $\leq\sqrt{2}$, and (\ref{ifbounded}) is satisfied with $M=2\sqrt{2}$.}; the coefficient was tuned empirically in pilot runs on small instances and is never changed afterwards. We also use two straightforward ``tricks'':
\begin{itemize}
\item Instead of considering one accuracy certificate, $\cC^N=\{y_t,\;\lambda^N_t=1/N,\;-\Psi(y_t)\}_{t=1}^N$, we build a ``bunch'' of certificates
\[
\cC_\mu^\nu=\bigg\{y_t,\;\lambda_t={1\over \nu-\mu+1},\;-\Psi(y_t)\bigg\}_{t=\mu}^\nu,
 \]where $\mu$ runs through a grid in $\{1,...,N\}$ (in this implementation, a 16-element equidistant grid), and $\nu\in\{\mu,\mu+1,...,N\}$ runs through another equidistant grid (e.g., for the largest problem instance, the grid $\{1,9, 17,...,257\}$). We compute the resolutions of these certificates and  identify the best (with the smallest resolution) certificate obtained so far.  Every 8 steps, the best certificate is used  to compute the current approximate solution to (\ref{saddlepoint}) along with the saddle point inaccuracy of this solution.
\item When applying MD to problem (\ref{saddlepoint}), the ``dual iterates'' $y_t=[\xi_t;\eta_t]$ and the ``primal iterates'' $x_t:=x(y_t)=[v_t;w_t]$ are pairs of matrices, with $n\times n$ matrices $\xi_t,\eta_t,v_t$ and $m\times m$ matrices $w_t$ (recall that we are in the case of $p_v=q_v=n$, $p_b=q_b=m$). It is easily seen that with $\cA$ given by (\ref{ourcA}), the matrices $\xi_t,\eta_t,v_t$ are linear combinations of rank 1 matrices $\alpha_i\beta_i^T$, $1\leq i\leq (k+1)t$, and $w_t$ are linear combinations of rank 1 matrices $\delta_i\epsilon_i^T$, $1\leq i\leq t$, with on-line computable vectors $\alpha_i,\beta_i,\delta_i,\epsilon_i$. Every step of MD adds $k+1$ new $\alpha$- and $k+1$ new $\beta$-vectors, and a pair of new $\delta$- and $\epsilon$-vectors. Our matrix iterates were represented by the vectors of coefficients in the above rank 1 decompositions (let us call this representation {\sl incremental}), so that the  computations performed at a step of MD, including computing the leading singular vectors by straightforward power iterations, are as if the standard representations of matrices were used, but all these matrices were of the size (at most) $n\times[(k+1)N]$, and not $n\times n$ and $m\times m$, as they actually are. In our experiments, for $k=2$ and $N\leq 512$, this incremental representation of iterates yields meaningful computational savings (e.g., by factor of $6$ for $n=8192$) as compared to the plain representation of iterates by 2D arrays.
\end{itemize}
\paragraph{Typical results} of our preliminary experiments are presented in Table \ref{table1}. There $\cC^t$ stands for the best certificate found in course of $t$ steps, and $\Gap(\cC^t)$ denotes the saddle point inaccuracy of the solution to (\ref{saddlepoint}) induced by this certificate (so that $\Gap(\cC^t)$ is a valid upper bound on the inaccuracy, in terms of the objective, to which the problem of interest (\ref{fineq1}) was solved in course of  $t$ steps). The comments are as follows:
    \begin{enumerate}
\item The results clearly demonstrate {\sl ``nearly linear''}, and not quadratic, growth of running time with $m,n$; this is due to the incremental  representation of iterates.
\item When evaluating the ``convergence patterns'' presented in the table, one should keep in mind that we are dealing with a method with slow $O(1/\sqrt{N})$ convergence rate, and from this perspective, 50-fold reduction in resolution in 512 steps is not that bad.
\item A natural alternative to the proposed approach would be to solve the saddle point problem (\ref{saddlepoint}) ``as it is,'' by applying to the associated primal v.i. (where the domain is the product of two nuclear norm balls and the operator is Lipschitz continuous and even skew symmetric) a proximal type saddle point algorithm and computing the required prox-mappings via full singular value decompositions. The state-of-the-art MP algorithm when applied to this problem exhibits $O(1/N)$ convergence rate;\footnote{For the primal v.i., (\ref{ifsmooth}) holds true for some $L>0$ and $M=0$. Moreover, with properly selected proximal setup for (\ref{fineq1}) the complexity bound (\ref{MPrate}) becomes $\Res(\cC^N|Y)\leq O(1)\sqrt{\ln(n)\ln(m)}/N$.} yet, every step of this method would require 2 SVD's of $n\times n$, and 2 SVD's of $m\times m$ matrices. As applied to the primal v.i.,
    MD exhibits $O(1/\sqrt{N})$ convergence rate, but the steps are cheaper -- we need one SVD of $n\times n$, and one SVD of an $m\times m$ matrix, and we are unaware of a proximal type algorithm for the primal v.i. with cheaper iterations. For the sizes $m,n,k$ we are interested in, the computational effort required by the outlined SVD's is, for all practical purposes, the same as the overall effort per step. Taking into account the actual SVD cpu times on the platform used in our experiments, the overall running times presented in Table \ref{table1}, i.e., times required by 512 steps of MD as applied to the dual v.i., allow for the following iteration counts $N$ for MP as applied to the primal v.i.:
    $$
    \begin{array}{c|c|c|c|c|}
    n&1024&2048&4096&8192\\
    \hline
    N&406&72&17&4\\
    \end{array}
    $$
    and for twice larger iteration counts for MD. From our experience, for $n=1024$ (and perhaps for $n=2048$ as well), MP algorithm as applied to the primal v.i. would yield solutions of better quality than those obtained with our approach. It, however, would hardly be the case, for both MP and MD, when $n=4096$, and definitely would not be the case for $n=8192$. Finally,
    with $n=16384$, CPU time used by the 257-step MD as applied to the dual v.i. is hardly enough to complete {\sl just one} iteration of MD as applied to the primal v.i. We believe these data
    demonstrate that the approach developed in this paper has certain practical potential.
\end{enumerate}
\begin{table}
\begin{center}
{\scriptsize
\begin{tabular}{||c|c||c|c|c|c|c|c|c|c|c||}
\cline{3-11}
\multicolumn{2}{c||}{}&\multicolumn{9}{|c||}{Iteration count $t$}\\
\cline{3-11}
\multicolumn{2}{c||}{}&1  &  65 &  129  & 193  & 257 &  321  & 385 &   449 &  512\\
\hline\hline
&$\Res(\cC^t|Y)$&1.5402 & 0.1535 & 0.0886 & 0.0621 & 0.0487 & 0.0389 & 0.03288 & 0.0293 & 0.0278\\
\cline{2-11}
$n=1024$&$\Res(\cC^1|Y)/\Res(\cC^t|Y)$&1.00 & 10.04 & 17.38 & 24.79 & 31.61 & 39.61 & 46.84 & 52.64 & 55.41\\
\cline{2-11}
$m=512$&$\Gap(\cC^t)$&0.1269 & 0.0239 & 0.0145 & 0.0103 & 0.0075 & 0.0063 & 0.0042 & 0.0040 & 0.0040\\
\cline{2-11}
$k=2$&$\Gap(\cC^t)/\Gap(\cC^t)$&1.00 & 5.31 & 8.78 & 12.38 & 17.03 & 20.20 & 29.98 & 31.41 & 31.66\\
\cline{2-11}
&cpu, sec&0.2 &      9.5 &     27.6 &     69.1 &    112.6 &    218.1 &    326.2 &    432.6 &    536.4\\
\hline\hline
&$\Res(\cC^t|Y)$&1.4809 & 0.1559 & 0.0842 & 0.0607 & 0.0471 & 0.0391 & 0.0337 & 0.0306 & 0.0285\\
\cline{2-11}
$n=2048$
&$\Res(\cC^1|Y)/\Res(\cC^t|Y)$&1.00 & 9.50 & 17.59 & 24.38 & 31.43 & 37.88 & 43.89 & 48.36 & 51.96\\
\cline{2-11}
$m=1024$&$\Gap(\cC^t)$&0.1329 & 0.0196 & 0.0119 & 0.0075 & 0.0053 & 0.0041 & 0.0036 & 0.0034 & 0.0027\\
\cline{2-11}
$k=2$&$\Gap(\cC^t)/\Gap(\cC^t)$&1.00 & 6.79 & 11.21 & 17.81 & 25.09 & 32.29 & 37.23 & 38.70 & 50.06\\
\cline{2-11}
&cpu, sec&0.7 &     38.0 &    101.1 &    206.3 &    314.1 &    508.9 &    699.0 &    884.9 &   1070.0\\
\hline\hline
&$\Res(\cC^t|Y)$&1.4845 & 0.1476 & 0.0891 & 0.0605 & 0.0491 & 0.0395 & 0.0329 & 0.0292 & 0.0275\\
\cline{2-11}
$n=4096$&$\Res(\cC^1|Y)/\Res(\cC^t|Y)$& 1.00 & 10.06 & 16.66 & 24.53 & 30.25 & 37.60 & 45.17 & 50.85 & 53.95\\
\cline{2-11}
$m=2048$&$\Gap(\cC^t)$& 0.1239 & 0.0222 & 0.0139 & 0.0108 & 0.0086 & 0.0041 & 0.0037 & 0.0035 & 0.0035\\
\cline{2-11}
$k=2$&$\Gap(\cC^t)/\Gap(\cC^t)$& 1.00 & 5.57 & 8.93 & 11.48 & 14.40 & 30.48 & 33.14 & 35.76 & 35.77\\
\cline{2-11}
&cpu, sec&  2.2 &    103.5 &    257.6 &    496.9 &    742.5 &   1147.8 &   1564.4 &   1981.4 &   2401.0\\
\hline\hline
&$\Res(\cC^t|Y)$& 1.4778 & 0.1391 & 0.0888& 0.0590 & 0.0469 & 0.0386 & 0.0324 & 0.0289 & 0.0270\\
\cline{2-11}
$n=8192$&$\Res(\cC^1|Y)/\Res(\cC^t|Y)$& 1.00 & 10.63 & 16.64 & 25.06 & 31.53 & 38.29 & 45.68 & 51.10 & 54.76\\
\cline{2-11}
$m=4096$&$\Gap(\cC^t)$& 0.1193 & 0.0232 & 0.0134 & 0.0108 & 0.0054 & 0.0040 & 0.0035 & 0.0034 & 0.0034\\
\cline{2-11}
$k=2$&$\Gap(\cC^t)/\Gap(\cC^t)$& 1.00 & 5.14 & 8.90 & 11.08 & 22.00 & 29.83 & 33.93 & 34.85 & 35.14\\
\cline{2-11}
&cpu, sec&      6.5 &    289.9 &    683.8 &   1238.1 &   1816.0 &   2724.5 &   3648.3 &   4572.2 &   5490.8\\
\hline\hline
&$\Res(\cC^t)$&1.4566 & 0.1154 & 0.0767 & 0.0556 & 0.0447&\multicolumn{4}{|c}{}\\
\cline{2-7}
$n=16384$&$\Res(\cC^1|Y)/\Res(\cC^t|Y)$& 1.00 & 12.62 & 19.00 & 26.22 & 32.60&\multicolumn{4}{|c}{}\\
\cline{2-7}
$m=8192$&$\Gap(\cC^t)$& 0.11959 & 0.02136 & 0.01460 & 0.01011 & 0.00853&\multicolumn{4}{|c}{}\\
\cline{2-7}
$k=2$&$\Gap(\cC^t)/\Gap(\cC^t)$&1.00 & 5.60 & 8.19 & 11.82 & 14.01&\multicolumn{4}{|c}{}\\
\cline{2-7}
&cpu, sec& 21.7 &    920.4 &   2050.2 &   3492.4 &   4902.2&\multicolumn{4}{|c}{}\\
\cline{1-7}
\end{tabular}}
\end{center}
\caption{\label{table1} MD on problem (\ref{fineq1}). Platform: 3.40 GHz i7-3770 desktop with 16 GB RAM, 64 bit Windows 7 OS.}
\end{table}
\subsubsection{Experiments with the MP-based scheme}
 In this section we briefly report on the results obtained with the modified MP-based scheme presented in section  \ref{secmodif}. Same as above, we use the  test problems and representation (\ref{thedata}) of the monotone operator of interest (with the only difference that now $Y=F$), and the Euclidean proximal setup. Using the Euclidean setup on $Y=F$ makes prox-mappings and functions $f_y(\cdot)$, defined in \rf{fyofx}, extremely simple:
$$
\begin{array}{rcl}
\Prox_{[\xi;\eta]}([d\xi;d\eta])&=&[\xi-d\xi;\eta-d\eta]\quad[\xi,d\xi,\eta,d\eta\in\bR^{p_v\times q_v}]\\
f_{y}(x)&=&{1\over 2}\langle y-\gamma[Gy-A^*x],y-\gamma[Gy-A^*x]\rangle+\gamma\langle a,x\rangle\\
&=&{1\over 2}\left[\|\xi+\gamma\eta+\gamma\cA^*w\|_\Fro^2+\|\eta-\gamma\xi+\gamma v\|_\Fro^2\right]+\gamma\langle b,w\rangle_\Fro\\
&&y=[\xi;\eta],x=[v;w].\\
\end{array}
$$
When choosing $\|\cdot\|_E$ to be the Frobenius norm,
\[
\|\underbrace{[v;w]}_x\|_E=\|x\|_{E,*}=\sqrt{\|v\|_\Fro^2+\|w\|_\Fro^2}.
\]
and taking into account that the spectral norm of $\cA$ is $\leq1$, it is immediately seen that the quantities $L$, $\gamma$, $\cL$ introduced in section \ref{sec:strategy}, can be set to 1, and what was called $R_E(X)$ in section \ref{sec:complexanal}, can be set to $\sqrt{2}$. As a result, by the complexity analysis of section \ref{sec:complexanal}, in order to find an $\epsilon$-solution to the problem of interest, we need $O(1)\epsilon^{-1}$ iterations of the recurrence (\ref{recurrence}), with $O(1)\epsilon^{-1}$ CGA steps of minimizing $f_{y_t}(\cdot)$ over $X$ per iteration, that is, the total of at most $O(1)\epsilon^{-2}$ calls to the LMO for $X=V\times W$. In fact, in our implementation $\epsilon$ is not fixed in advance; instead, we fix the total number $N=256$ of calls to LMO, and terminate CGA at iteration $t$ of the recurrence (\ref{recurrence}) when either a solution $x_t\in x$ with $\delta_t\leq 0.1/t$ is achieved, or the number of CGA steps reaches a prescribed limit (set to 32 in the experiment to be reported).
\par
Same as in the first series of experiments, ``incremental'' representation of matrix iterates is used in the experiments with the MP-based scheme. In these experiments we also use a special {\sl post-processing} of the solution we explain next.
\paragraph{Post-processing.} Recall that in the situation in question the step $\# i$ of the CGA at iteration $\# t$ of the MP-based recurrence produces a pair $[v_{t,i};w_{t,i}]$ of rank 1 of $n\times n$ and $m\times m$ matrices  of unit spectral norm -- the minimizers of the linear  form $\langle \nabla f_{y_t}(x_{t,i}),x\rangle$ over $x\in X$; here $x_{t,i}$ is $i$-th step of CGA minimization of $f_{y_t}(\cdot)$ over $X$. As a result, upon termination, we have at our disposal $N=256$ pairs of rank one matrices $[v_j;w_j]$, $1\leq j\leq N$, known to belong to $X$. Note that the approximate solution $\widehat{x}$, as defined in section \ref{sec:constr5}, is a certain convex combination of these matrices. A natural way to get a better solution is to solve the optimization problem
\be
\Opt=\min_{\lambda,v}\left\{f(\lambda)=\|\cA v-b\|_{2,2}: v={\sum}_{j=1}^N\lambda_j v_j,{\sum}_{j=1}^N|\lambda_j|\leq1\right\}.
\ee{(Q)}
Indeed, note that the $v$-components of feasible solutions to this problem are of nuclear norm $\leq1$, i.e., are feasible solutions to the problem of interest (\ref{fineq1}), and that in terms of the objective of (\ref{fineq1}), the $v$-component of an optimal solution to \rf{(Q)} can be only better than the $v$-component of $\widehat{x}$. On the other hand, \rf{(Q)} is a low-dimensional convex optimization problem on a simple domain, and the first order information on $f$ can be obtained, at a relatively low cost, by Power Method, so that \rf{(Q)} is well suited for solving by proximal first order algorithms, e.g., the Bundle Level algorithm \cite{LNN} we use in our experiments.
\paragraph{Numerical illustration.} Here we present just one (in fact, quite representative) numerical example. In this example $n=4096$ and $m=2048$ (i.e., in (\ref{fineq1}) the variable matrix $u$  is of size $4096\times 4096$, and the data matrix $b$ is of size $2048\times 2048$); the mapping $\cA$ is given by (\ref{ourcA}) with $k=2$. The data are generated in the same way as in the experiments described in section \ref{sec:exp-MD} except for the fact that we used $b=\cA \bar{u}$ to ensure zero optimal value in (\ref{fineq1}). As a result, the value of the objective of (\ref{fineq1}) at an  approximate solution coincides with the inaccuracy of this solution in terms of the objective of (\ref{fineq1}). In the experiment we report on here, the objective of (\ref{fineq1}) evaluated at the initial -- zero -- solution, i.e., $\|b\|_{2,2}$, is equal to 0.751. After the total of 256 calls to the LMO for $X$ (just 11 steps of recurrence (\ref{recurrence})) and post-processing which took 24\% of the overall CPU time, the value of the objective is reduced to  0.013 -- by factor 57.3. For comparison, when processing the same instance by the basic MD scheme, augmented by the just outlined post-processing, after 256 MD iterations (i.e., after the same as above 256 calls to the LMO), the value of the objective at the resulting feasible solution to (\ref{fineq1}) was 0.071, meaning the progress in accuracy by factor 10.6 (5 times worse than the progress in accuracy for the MP-based scheme). Keeping the instance intact and increasing the number of MD iterations in the basic scheme from 256 to 512, the objective at the approximate solution yielded by the post-processing reduces from 0.071 to 0.047, which still is 3.6 times worse than that achieved with the MP-based scheme after 256 calls to LMO.

\appendix\section{Proofs}
\subsection{Proof of Theorems \ref{propmain} and \ref{propmainnew}}\label{sec:ensure}
We start with proving Theorem \ref{propmainnew}. In the notation of the theorem, we have
\be
\begin{array}{lcl}
&&\forall x\in X: \Phi(x)=Ay(x)+a,\\
(a):&&y(x)\in Y,\\
(b): &&\langle y(x)-y,A^*x-G(y(x))\rangle \geq 0\,\forall y\in Y.
\end{array}
\ee{proof-def}
For $\bar{x}\in X$, let $\bar{y}=y(\bar{x})$, and let $\widehat{y}=\sum_t\lambda_ty_t$, so  that $\bar{y},\widehat{y}\in Y$ by (\ref{proof-def}.a). Since $G$ is monotone, for all $t\in \{1,...,N\}$ we have
$$
\begin{array}{lll}
&&\langle \bar{y}-y_t,G(\bar{y})-G(y_t)\rangle \geq 0\\
\Rightarrow&&\langle \bar{y},G(\bar{y})\rangle \geq \langle y_t,G(\bar{y})\rangle +\langle \bar{y},G(y_t)\rangle -\langle y_t,G(y_t)\rangle\,\,\forall t\\
\Rightarrow& &\langle \bar{y},G(\bar{y})\rangle \geq\sum_t\lambda_t\left[\langle y_t,G(\bar{y})\rangle +\langle \bar{y},G(y_t)\rangle -\langle y_t,G(y_t)\rangle\right]\\
&&\hbox{[since $\lambda_t\geq0$ and $\sum_t\lambda_t=1$]},
\end{array}
$$
and we conclude that
\be
\langle \bar{y},G(\bar{y})\rangle-\langle \widehat{y},G(\bar{y})\rangle\geq \sum_{t=1}^N\lambda_t\left[\langle \bar{y},G(y_t)\rangle -\langle y_t,G(y_t)\rangle\right].
\ee{ee}
We now have
$$
\begin{array}{ll}
&\langle\Phi(\bar{x}),\bar{x}-\sum_t\lambda_tx_t\rangle\\
=&\langle A\bar{y}+a,\bar{x}-\sum_t\lambda_tx_t\rangle
=\langle \bar{y},A^*\bar{x}-\sum_t\lambda_tA^*x_t\rangle +\langle a,\bar{x}-\sum_t\lambda_tx_t\rangle\\
=&\langle \bar{y},A^*\bar{x}-G(\bar{y})\rangle +\langle \bar{y},G(\bar{y})-\sum_t\lambda_tA^*x_t\rangle+\langle a,\bar{x}-\sum_t\lambda_tx_t\rangle\\
\geq& \langle \widehat{y},A^*\bar{x}-G(\bar{y})\rangle +\langle \bar{y},G(\bar{y})-\sum_t\lambda_tA^*x_t\rangle+\langle a,\bar{x}-\sum_t\lambda_tx_t\rangle\\
&\hbox{\ [by (\ref{proof-def}.b) with $y=\widehat{y}$ and due to $\bar{y}=y(\bar{x})$]}\\
=&\langle \widehat{y},A^*\bar{x}\rangle +\left[\langle G(\bar{y}),\bar{y}\rangle-\langle G(\bar{y}),\widehat{y}\rangle\right]-\langle \bar{y},\sum_t\lambda_tA^*x_t\rangle
+\langle a,\bar{x}-\sum_t\lambda_tx_t\rangle\\
\geq&\langle \widehat{y},A^*\bar{x}\rangle +\sum_t\lambda_t\left[\langle \bar{y},G(y_t)\rangle -\langle y_t,G(y_t)\rangle\right]-\langle \bar{y},\sum_t\lambda_tA^*x_t\rangle
+\langle a,\bar{x}-\sum_t\lambda_tx_t\rangle\,\,\hbox{[by \rf{ee}]}\\
=&\sum_t\lambda_t\langle y_t,A^*\bar{x}\rangle+\sum_t\lambda_t\left[\langle \bar{y},G(y_t)\rangle -\langle y_t,G(y_t)\rangle-\langle \bar{y},A^*x_t\rangle
+\langle a,\bar{x}-x_t\rangle\right]\\
&\hbox{[since $\widehat{y}=\sum_t\lambda_ty_t$ and $\sum_t\lambda_t=1$]}\\
=&\sum_t\lambda_t\left[\langle Ay_t,\bar{x}-x_t\rangle+\langle Ay_t,x_t\rangle
+\langle \bar{y},G(y_t)\rangle -\langle y_t,G(y_t)\rangle-\langle \bar{y},A^*x_t\rangle
+\langle a,\bar{x}-x_t\rangle\right]\\
=&\sum_t\lambda_t\left[\langle y_t,A^*x_t\rangle +\langle \bar{y},G(y_t)\rangle -\langle y_t,G(y_t)\rangle-\langle \bar{y},A^*x_t\rangle\right]+
\sum_t\lambda_t\langle Ay_t+a,\bar{x}-x_t\rangle\\
=&{\sum}_t\lambda_t\langle A^*x_t-G(y_t),y_t-\bar{y}\rangle
+{\sum}_t\lambda_t\langle Ay_t+a,\bar{x}-x_t\rangle
\ge -\epsilon+{\sum}_t\lambda_t\langle Ay_t+a,\bar{x}-x_t\rangle\\
&~~~~~~~~~~~~~~~~~~~~~~~~~~~~~~~~~~~~~~~~~~~~~~~~~~~~~~~~~~~~~~~\hbox{[by 
\rf{epsilonetc} due to $\bar{y}=y(\bar{x})\in Y(X)$].}
\end{array}
$$
The bottom line is that
$$
\langle \Phi(\bar{x}),\widehat{x}-\bar{x}\rangle \leq\epsilon+{\sum}_{t=1}^N\lambda_t\langle Ay_t+a,x_t-\bar{x}\rangle\,\forall \bar{x}\in X,
$$
as stated in (\ref{bottomline}). Theorem \ref{propmainnew} is proved.\par
To prove Theorem \ref{propmain}, let $y_t\in Y$, $1\leq t\leq N$, and $\lambda_1,...,\lambda_N$ be from the premise of the theorem, and let $x_t$, $1\leq t\leq N$, be specified as
$x_t=x(y_t)$, so that $x_t$ is the minimizer of the linear form $\langle Ay_t+a,x\rangle$ over $x\in X$. Due to the latter choice, we have $\sum_{t=1}^N\lambda_t\langle Ay_t+a,x_t-\bar{x}\rangle\leq0$ for all $\bar{x}\in X$, while $\epsilon$ as defined by (\ref{epsilonetc}) is nothing but $\Res(\{y_t,\lambda_t,-\Psi(x_t)\}_{t=1}^N|Y(X))$. Thus, (\ref{bottomline}) in the case in question implies that
$$
\forall \bar{x}\in X: \langle \Phi(\bar{x}),{\sum}_{t=1}^N\lambda_t x_t - \bar{x}\rangle \leq \Res(\{y_t,\lambda_t,-\Psi(x_t)\}_{t=1}^N|Y(X)),
$$
and (\ref{ensure}) follows. Relation (\ref{ensureskew}) is an immediate corollary of (\ref{ensure}) and Lemma \ref{lem2} as applied to $X$ in the role of $Y$, $\Phi$  in the role of $H$, and $\{x_t,\lambda_t,\Phi(x_t)\}_{t=1}^N$ in the role of $\cC^N$. \qed
\subsection{Proof of Proposition \ref{propMD}}
Observe that the optimality conditions in the optimization problem specifying $v=\Prox_y(\zeta)$ imply that
\[
\langle \xi-\omega'(y)+\omega'(v),z-v\rangle \geq0,\,\,\forall z\in Y,
\]or
\[
\langle \xi,v-z\rangle\leq \langle \omega'(v)-\omega'(y),z-v\rangle=\langle V'_y(v),z-v\rangle,\,\,\forall z\in Y,
\]
which, using a remarkable identity \cite{CT93}
\[
\langle V'_y(v),z-v\rangle=V_y(z)-V_v(z)-V_y(v),
\]
can be rewritten equivalently as
\begin{equation}\label{app1}
 v=\Prox_y(\zeta)\Rightarrow \langle\zeta,v-z\rangle \leq V_y(z)-V_v(z)-V_y(v)\,\,\forall z\in Y.
\end{equation}
Setting $y=y_t$, $\xi=\gamma_t H_t(y_t)$, which results in $v=y_{t+1}$, we get
$$
\forall z\in Y: \gamma_t\langle H_t(y_t),y_{t+1}-z\rangle\leq V_{y_t}(z)-V_{y_{t+1}}(z)-V_{y_t}(y_{t+1}),
$$
whence,
\bse
\forall z\in Y: \gamma_t\langle H_t(y_t),y_{t}-z\rangle&\leq& V_{y_t}(z)-V_{y_{t+1}}(z)+\underbrace{\left[\gamma_t\langle H_t(y_t),y_t-y_{t+1}\rangle-V_{y_t}(y_{t+1})\right]}_{\leq \gamma_t\|H_t(y_t)\|_*\|y_t-y_{t+1}\|-\half\|y_t-y_{t+1}\|^2}\\
&\leq& V_{y_t}(z)-V_{y_{t+1}}(z)+\half\gamma_t^2\|H_t(y_t)\|_*^2.
\ese
Summing up these inequalities over $t=1,...,N$ and taking into account that for $z\in Y'$, we have $V_{y_1}(z)\leq \half\Omega^2[Y']$ and that $V_{y_{N+1}}(z)\geq0$, we get (\ref{thenMD}). \qed
\subsection{Proof of Proposition \ref{propMP}}
Applying (\ref{app1}) to $y=y_t$, $\xi=\gamma_t H_t(z_t)$, which results in $v=y_{t+1}$, we get
$$
\forall z\in Y: \gamma_t\langle H_t(z_t),y_{t+1}-z\rangle \leq V_{y_t}(z)-V_{y_{t+1}}(z)-V_{y_t}(y_{t+1}),
$$
whence, by the definition \rf{dt} of $d_t$,
\begin{equation}\label{verylast-1}
\begin{array}{l}
\forall z\in Y: \gamma_t\langle H_t(z_t),z_t-z\rangle\leq V_{y_t}(z)-V_{y_{t+1}}(z)+d_t.\\
\end{array}
\end{equation}
Summing up the resulting inequalities over $t=1,...,N$ and taking into account that $V_{y_1}(z)\leq \half\Omega^2[Y']$ for all $z\in Y'$ and
$V_{y_{N+1}}(z)\geq0$, we get
$$
\forall z\in Y': \sum_{t=1}^n\lambda^N_t\langle H_t(z_t),z_t-z\rangle \leq {\half\Omega^2[Y']+\sum_{t=1}^N d_t\over\sum_{t=1}^N\gamma_t}.
$$
The right hand side in the latter inequality is independent of $z\in Y'$. Taking supremum of the left hand side over $z\in Y'$, we arrive at
\rf{thenMP1}.
\par
Moreover, invoking (\ref{app1}) with $y=y_t$, $\xi=\gamma_t H_t(y_t)$ and specifying $z$ as $y_{t+1}$, we get
$$
\gamma_t\langle H_t(y_t),z_t-y_{t+1}\rangle \leq V_{y_t}(y_{t+1})-V_{z_t}(y_{t+1})-V_{y_t}(z_t),
$$
whence
\begin{equation}\label{verylast}
\begin{array}{rcl}
d_t&=&\gamma_t\langle H_t(z_t),z_t-y_{t+1}\rangle - V_{y_t}(y_{t+1})\leq \gamma_t\langle H_t(y_t),z_t-y_{t+1}\rangle
+\gamma_t\langle H_t(z_t)-H_t(y_t),z_t-y_{t+1}\rangle\\
&&- V_{y_t}(y_{t+1})\\
&\leq& -V_{z_t}(y_{t+1})-V_{y_t}(z_t)+\gamma_t\langle H_t(z_t)-H_t(y_t),z_t-y_{t+1}\rangle\\
&\leq&
\gamma_t\|H_t(z_t)-H_t(y_t)\|_*\|z_t-y_{t+1}\|-{\half}\|z_t-y_{t+1}\|^2-{\half}\|y_t-z_t\|^2\\
&\leq& {\half}\left[\gamma_t^2\|H_t(z_t)-H_t(y_t)\|_*^2-\|y_t-z_t\|^2\right],\\
\end{array}
\end{equation}
as required in \rf{thenMP2}. \qed
\subsection{Proof of Lemma \ref{lemmaholdstrue}}\label{Proofoflemmaholdstrue}
{\bf 1$^0$.} We start with the following standard fact:
\begin{lemma}\label{wellknownlemma} Let $Y$ be a nonempty closed convex set in Euclidean space $F$, $\|\cdot\|$ be a norm on $F$, and $\omega(\cdot)$ be a continuously differentiable function on $Y$ which is strongly convex, modulus 1, w.r.t. $\|\cdot\|$. Given $b\in F$ and $y\in Y$, let us set
$$
\begin{array}{rcl}
g_y(\xi)&=&\max\limits_{z\in Y}\left[\langle z,\omega'(y)-\xi\rangle -\omega(z)\right]:F\to \bR,\\
z_y(\xi)&=&\argmax\limits_{z\in Y}\left[\langle z,\omega'(y)-\xi\rangle -\omega(z)\right].\\
\end{array}
$$
The function $g_y$ is convex with Lipschitz continuous gradient $\nabla g_y(\xi)=-z_y(\xi)$:
\begin{equation}\label{Lipschitzcont}
\|\nabla g_y(\xi)-\nabla g_y(\xi')\|\leq \|\xi-\xi'\|_*\,\,\forall \xi,\xi',
\end{equation}
where $\|\cdot\|_*$ is the norm conjugate to $\|\cdot\|$.
\end{lemma}
Indeed, since $\omega$ is strongly convex and continuously differentiable on $Y$, $z_y(\cdot)$ is well defined, and from optimality conditions it holds
\begin{equation}\label{eq987654}
\langle \omega'(z_y(\xi))+\xi-\omega'(y),z_y(\xi)-z\rangle \leq0\,\,\forall z\in Y.
\end{equation}
Consequently, $g_y(\cdot)$  is well defined; this function clearly is convex, and the vector $-z_y(\xi)$ clearly is a subgradient of $g_y$ at $\xi$. If now $\xi',\xi''\in F$, then, setting $z'=z_y(\xi')$, $z''=z_y(\xi'')$ and invoking (\ref{eq987654}), we get
$$
\langle \omega'(z')+\xi'-\omega'(y),z'-z''\rangle \leq0,\,\,\langle\omega'(z'')+\xi''-\omega'(y),z''-z'\rangle \leq0
$$
whence, summing the inequalities up,
$$
\langle \xi'-\xi'',z'-z''\rangle \leq \langle \omega'(z')-\omega'(z''),z''-z'\rangle \leq -\|z'-z''\|^2,
$$
implying that $\|z'-z''\|\leq \|\xi'-\xi''\|_*$. Thus, a subgradient field $-z_y(\cdot)$ of $g_y(\cdot)$ is Lipschitz continuous with constant 1 from $\|\cdot\|_*$ into $\|\cdot\|$, whence $g_y$ is continuously differentiable and (\ref{Lipschitzcont}) takes place. \qed
\indent{\bf 2$^0$.} To derive Lemma \ref{lemmaholdstrue} from Lemma \ref{wellknownlemma}, set in the latter Lemma $Y=F$ and note that
 $f_y(x)$ is obtained from $g_y(\cdot)$ by affine substitution of variables and adding linear form:
$$
f_y(x)=g_y(\gamma[Gy- A^*x])+\gamma\langle a,x\rangle.
$$
whence $\nabla f_y(x)=-\gamma A\nabla g_y(\gamma[Gy-A^*x])+\gamma a=\gamma A z_y(\gamma[Gy-A^*x])+\gamma a$, as required in (\ref{eqgradf}), and
$$
\begin{array}{rcl}
\|\nabla f_y(x')-\nabla f_y(x'')\|_{E,*}&=&\gamma \|A[z_y(Gy-A^*x')-z_y(Gy-A^*x'')]\|_{E,*}\\
&\leq& (\gamma\cL)\|\nabla g_y(\gamma[Gy-A^*x'])-\nabla g_y(\gamma[Gy-A^*x''])\|\\
&\leq&  (\gamma\cL)\|\gamma[Gy-A^*x']-\gamma[Gy-A^*x'']\|_*\leq (\gamma \cL)^2\|x'-x''\|_E\\
\end{array}
$$
(we have used (\ref{Lipschitzcont}) and equivalences in (\ref{calLnorm})), as required in (\ref{LipContGrad}). \qed

\subsection{Review of Conditional Gradient Algorithm}\label{sec:CGA}
  The required description of CGA and its complexity analysis are as follows.\par
  As applied to minimizing a smooth -- with Lipschitz continuous gradient
$$
\|\nabla f(u)-\nabla f(u')\|_{E,*}\leq \cL\|u-u'\|_E,\,\,\forall u,u'\in X,
$$
convex function $f$ over a convex compact set $X\subset E$, the generic CGA is the recurrence of the form
$$
\begin{array}{rcl}
u_1&\in& X\\
u_{s+1}&\in& X \hbox{\ satisfies\ } f(u_{s+1})\leq f(u_s+\gamma_s [u_s^+-u_s]),\,s=1,2,...\\
&&\gamma_s={2\over s+1},\,u_s^+\in\Argmin_{u\in X}\langle f'(u_s),u\rangle.\\
\end{array}
$$
The standard results on this recurrence (see, e.g., proof of Theorem 1 in \cite{HJN}) state that if $f_*=\min_X f$, then
\begin{equation}\label{steq123}
\begin{array}{ll}
(a)&\epsilon_{t+1}:=f(u_{t+1})-f_*\leq\epsilon_{t}-\gamma_t\delta_t+{2\cL R^2\gamma_t^2},\,t=1,2,...\\
&\delta_t:=\max_{u\in X} \langle \nabla f(u_t),u_t-u\rangle;\\
(b)&\epsilon_t\leq {2\cL R^2\over t+1},t=2,3,...\\
\end{array}
\end{equation}
where $R$ is the smallest of the radii of $\|\cdot\|_E$-balls containing $X$.
From (\ref{steq123}.$a$) it follows that
$$
\gamma_\tau\delta_\tau\leq \epsilon_\tau-\epsilon_{\tau+1}+2\cL R^2\gamma_\tau^2,\,\tau=1,2,...;
$$
summing up these inequalities over $\tau=t,t+1,...,2t$, where $t>1$, we get
$$
\left[\min_{\tau\leq 2t}\delta_\tau\right]\sum_{\tau=t}^{2t} \gamma_\tau\leq \epsilon_t+2\cL R^2\sum_{\tau=t}^{2t}\gamma_\tau^2,
$$
which combines with (\ref{steq123}.$b$) to imply that
$$
\min_{\tau\leq 2t}\delta_\tau \leq O(1)\cL R^2{{1\over t}+\sum_{\tau=t}^{2t}{1\over \tau^2}\over\sum_{\tau=t}^{2t}{1\over\tau}}\leq O(1){\cL R^2\over t}.
$$
It follows that given $\epsilon<\cL R^2$, it takes at most $O(1){\cL R^2\over \epsilon}$ steps of CGA to generate a point $u^\epsilon\in X$
with $\max_{u\in X} \langle \nabla f(u^\epsilon),u^\epsilon-u\rangle \leq\epsilon$.
\end{document}